\documentclass[11pt,letterpaper,oneside,leqno]{article}

\usepackage[T1]{fontenc}
\usepackage[utf8]{inputenc}

\usepackage[letterpaper, top=1in, bottom=1.1in]{geometry}
\usepackage[numbers]{natbib}
\usepackage{booktabs}
\usepackage[all]{nowidow}
\usepackage[labelfont={bf,small}, labelsep=period]{caption}

\usepackage{amsthm,amsmath,amssymb}
\usepackage{mathtools}
\usepackage{bm}
\usepackage{siunitx}

\usepackage{graphicx}
\usepackage[colorlinks,citecolor=blue,urlcolor=blue,linkcolor=blue,linktocpage=true]{hyperref}

\numberwithin{equation}{section}
\theoremstyle{plain}

\newtheorem{theorem}{Theorem}
\newtheorem{lemma}{Lemma}
\newtheorem{corollary}{Corollary}
\newtheorem{definition}{Definition}

\newcommand{\z}{\bm{z}}
\newcommand{\p}{\bm{p}}
\newcommand{\q}{\bm{q}}

\newcommand{\Filter}{\mathcal{F}}

\newcommand{\Pbb}{\mathbb{P}}
\newcommand{\Qbb}{\mathbb{Q}}  
\newcommand{\Pcal}{\mathcal{P}}
\newcommand{\Qcal}{\mathcal{Q}}

\newcommand{\Rnng}{\mathbb{R}_{\geq 0}}

\newcommand{\PiLambda}{\Pi_{\bm{\lambda}}}
\newcommand{\PiAlpha}{\Pi_{\overline{\bm{\alpha}}}}

\newcommand{\Gfunc}{G_{\bm{\lambda}}^{(k)}}
\newcommand{\gfunc}{g_{\bm{\lambda}}}
\newcommand{\dfunc}{d_{\bm{\lambda}}}
\newcommand{\Dfunc}{D_{\bm{\lambda}}}

\newcommand{\rhofunc}{\rho_{\bm{\lambda}}}
\newcommand{\rhodomain}{\Omega_\rho}
\newcommand{\rhofilter}{\mathcal{F}_\rho}

\newcommand{\dint}{\mathrm{d}}

\newcommand{\titleref}[1]{\texorpdfstring{\MakeLowercase{\ref{#1}}}{\ref{#1}}}
\newcommand{\refeqq}[1]{\mathrel{\overset{\makebox[0pt]{\mbox{\normalfont\scriptsize\sffamily #1}}}{=}}}

\DeclareMathOperator*{\argmax}{arg\,max}

\let\originalleft\left
\let\originalright\right
\renewcommand{\left}{\mathopen{}\mathclose\bgroup\originalleft}
\renewcommand{\right}{\aftergroup\egroup\originalright}

\newcommand{\RN}[1]{\textup{\uppercase\expandafter{\romannumeral#1}}}

\title{Minimax Optimal Sequential Hypothesis Tests \\ for Markov Processes}
\author{
  Michael Fau\ss{}, H.~Vincent Poor\\
	Dept.~of Electrical Engineering\\
	Princeton University\\
	Princeton, NJ, USA \\
	\and 
  Abdelhak M.~Zoubir\\
	Signal Processing Group\\
	TU Darmstadt\\
	Darmstadt, Germany\\
}

\date{Updated October 2020}

\begin{document}

\maketitle

\begin{abstract}
  Under mild Markov assumptions, sufficient conditions for strict minimax optimality of sequential tests for multiple hypotheses under distributional uncertainty are derived. Two objective functions are considered, namely, the weighted sum of the expected run-length and the error probabilities, and the expected run-length under constraints on the error probabilities. First, the design of optimal sequential tests for multiple simple hypotheses is revisited and it is shown that the partial derivatives of the corresponding cost function are closely related to the performance metrics of the underlying sequential test. Second, an implicit characterization of the least favorable distributions for a given testing policy is stated. By combining the results on optimal sequential tests and least favorable distributions, sufficient conditions for a sequential test to be minimax optimal under general distributional uncertainties are obtained. The cost function of the minimax optimal test is further identified as a generalized $f$-dissimilarity and the least favorable distributions as those that are most similar with respect to this dissimilarity. Numerical examples for minimax optimal sequential tests under different distributional uncertainties illustrate the theoretical results.
\end{abstract}

\section{Introduction}
\label{sec:introduction}

Sequential hypothesis tests are well-known for being highly efficient in terms of the number of required samples and, as a consequence, for minimizing the decision delay in time-critical applications. In his seminal book \cite{Wald1947}, Wald showed that, compared to fixed sample size tests, sequential tests can reduce the average number of samples by a factor of two. In general, the ability to allow the overall number of samples to depend on the current history makes sequential procedures more flexible and adaptable than procedures whose sample size is chosen \emph{a priori}. Comprehensive overviews of sequential hypothesis testing and related topics can be found in \cite{Wald1947, Siegmund1985, Ghosh1991, Poor2009, Tartakovsky2014}, to name just a few.

A well-established drawback of sequential hypothesis tests is that their higher efficiency depends critically on the assumption that the process generating the observations indeed follows the assumed model. If this is not the case, that is, if a model mismatch occurs, the number of samples can increase significantly; compare, for example, \cite[Figure 3.4]{Tartakovsky2014}, which illustrates the influence of a model mismatch on the expected run length of a sequential test for the mean of an auto-regressive process. This observation lead Kiefer and Weiss to propose a sequential test that, in addition to meeting the targeted error probabilities under the hypotheses, minimizes the \emph{maximum} expected run length over all feasible distributions \cite{Dvoretzky1953, Kiefer1957}. Different variations of the corresponding optimization problem are known as Kiefer--Weiss problem or modified Kiefer--Weiss problem and have received considerable attention in the literature \cite{Lorden1976, Dragalin1988, Pavlov1991, Zhitlukhin2013}. However, to the present day, exact solutions to the Kiefer--Weiss problem have only been shown for special cases of binary hypothesis tests.

A natural generalization of the Kiefer--Weiss problem is to include the error probabilities in the minimax criterion, that is, to design a test whose maximum error probabilities are minimal over the set of feasible distributions. For fixed sample sizes this minimax approach to the design of statistical tests was pioneered by Huber \cite{Huber1965} and is known as \emph{robust hypothesis testing} or \emph{robust detection}. In general, robust hypothesis tests sacrifice some efficiency under ideal conditions in order to be less sensitive to deviations from the ideal case \cite{Huber1981}. In this sense, robust hypothesis tests, and robust statistics in general, form a middle ground between parametric and non-parametric approaches. For an overview of existing results, recent advances and applications of robust statistics, see, for example, \cite{KassamPoor1985, Maronna2006, Zoubir2012, Medina2015, Zoubir2018}.

The idea underlying this paper is to leverage both sequential and robust hypothesis testing. Ideally, a robust sequential test is \emph{fast and reliable}, that is, it requires fewer observations, on average, than a fixed sample size test and at the same time works reliably under model mismatch. The idea underlying this paper is to leverage both sequential and robust hypothesis testing. Ideally, a robust sequential test is \emph{fast and reliable}, i.e., it requires fewer observations, on average, than a fixed sample size test and at the same time works reliably under model mismatch. In what follows, a coherent framework for the design of minimax optimal sequential tests is developed, including the Kiefer--Weiss problem as a special case. The main contributions are summarized as follows:
\begin{itemize}
  \item Sufficient conditions for strict minimax optimality of sequential tests for multiple hypotheses under distributional uncertainty are derived. The results hold for stochastic processes with a time-homogeneous Markovian representation and under very mild assumptions on the type of uncertainty. 
  \item The tests satisfying the minimax optimality conditions are identified as optimal tests for least favorable distributions. The latter are shown to be data dependent such that the underlying stochastic process becomes Markovian, even if the increments are independent under the nominal model.
  \item It is shown that the least favorable distributions are those that are maximally similar with respect to a weighted statistical similarity measure of the $f$-dissimilarity type with the weights being data dependent. Both the similarity measure and the way in which the weights depend on the data are characterized explicitly.
  \item A method for the design of minimax optimal sequential tests in practice is proposed, and two numerical examples are provided to illustrate the theoretical findings.
\end{itemize}

The existing literature on minimax hypothesis testing can roughly be divided into two groups. The first group consists of results on strictly minimax optimal tests that are based on 2-alternating Choquet capacities \cite{Choquet1954, Huber1973}. If the existence of capacity achieving distributions can be shown, then these distributions are least favorable and maximally similar with respect to \emph{all} convex divergences; compare \cite[Theorem 6.1]{Huber1973}. The derivation of strictly minimax fixed sample size tests in \cite{Huber1965, Kassam1981, Fauss2016_old_bands} as well as the strictly minimax change detection procedure in \cite{Unnikrishnan2011} are based on this result.

If no capacity achieving distributions can be identified, least favorable distributions are usually defined as the minimizers of a given statistical divergence, playing the role of a surrogate objective function. That is, instead of solving the robust detection problem exactly, it is shown that choosing distributions that are maximally similar with respect to a suitable divergence are \emph{almost} least favorable and the corresponding tests \emph{almost} minimax. Typical examples are \cite{Poor1980, Basu1998, Pardo2005, Yilmaz2016, Gao2018}. However, many asymptotic results can also be grouped under this approach in the sense that the solutions are obtained by solving an asymptotic surrogate problem which typically induces the Kullback--Leibler divergence as the corresponding similarity measure. Examples are abundant in the literature: \cite{Degroot1960, Schmitz1987, Dragalin1988, El-Sawy1979, Voudouri1988, Brodsky2008, Fellouris2012, Zhitlukhin2013}. Comparable asymptotic result can also be found for the closely related problem of minimax quickest change detection \cite{Brodsky2008a, Fellouris2012a, Banerjee2015}; the connection between the two probelms is discussed in more detail in Section~\ref{sec:discussion}.

The few works on minimax sequential testing that do not fall under the two categories of asymptotic results or capacity based results are usually limited in scope. For example, the results in \cite{Vos2001} are highly application specific, the results in \cite{Maurice1957, Kharin2002} are based on strong assumptions on the distributions and the results in \cite{Gul2016} only hold for the single sample case. 

In summary, what sets this paper apart from the existing literature is that \emph{strictly minimax optimal} sequential tests are derived under \emph{mild distributional assumptions} and that the proof of optimality is given \emph{without invoking Choquet capacities}. Finally, note that the paper extends and generalizes \cite{Fauss2016_thesis} where minimax optimal sequential tests for two hypotheses are studied under slightly stricter assumptions.

The paper is organized as follows. The minimax sequential testing problem is stated formally in Section~\ref{sec:problem_formulation}. Optimal sequential tests for multiple simple hypotheses are revisited in Section~\ref{sec:optimal_tests}. Some useful results that relate properties of an optimal test to the partial derivatives of its cost function are shown in Section~\ref{sec:properties_of_rho}. Least favorable distributions and their properties are discussed in Section~\ref{sec:least_favorable_distributions}. The main result, sufficient conditions for minimax optimal sequential tests, is stated in Section~\ref{sec:minimax_optimal_sequential_tests} followed by a discussion in Section~\ref{sec:discussion}. Two illustrative numerical examples are shown in Section~\ref{sec:examples}.

\section{Notation and Problem Formulation}
\label{sec:problem_formulation}

In this section the minimax optimal sequential testing problem is defined, and some common notations are introduced. Notations not covered here are defined when they occur in the text.

\subsection{Notation}
\label{ssec:notation}

Random variables are denoted by uppercase letters and their realizations by lowercase letters. Analogously, probability distributions are denoted by uppercase letters and their densities by the corresponding lowercase letters. Blackboard bold is used to indicate product measures. Measurable sets are denoted by tuples $(\Omega, \Filter)$. Boldface lowercase letters are used to indicate vectors; no distinction is made between row and column vectors. The inner product of two vectors, $\bm{x}$ and $\bm{y}$, is denoted by $\langle \bm{x}, \bm{y} \rangle$ and the element-wise product by $\bm{xy}$. All comparisons between vectors are defined element-wise. The indicator function of a set $\mathcal{A}$ is denoted by $\mathcal{I}(\mathcal{A})$. All comparisons between functions are defined point-wise. 

The notation $\partial_{y_k} f(\bm{y})$ is used for the subdifferential \cite[\S 23]{Rockafellar1970} of a convex function $f \colon \mathcal{Y} \subset \mathbb{R}^K \to \mathbb{R}$ with respect to $y_k$ evaluated at $\bm{y}$, that is,
\begin{equation}
  \partial_{y_k} f(\bm{y}) \coloneqq \left\{\, c \in \mathbb{R} : f(\bm{y}') - f(\bm{y}) \geq c \, (y'_k-y_k) \quad \forall \bm{y}' \in \mathcal{Y}_k(\bm{y}) \, \right\},
  \label{eq:partial_diff}
\end{equation}
with $\mathcal{Y}_k(\bm{y}) \coloneqq \left\{\, \bm{y}' \in \mathcal{Y} : y_j' = y_j \quad \forall j \in \{1, \ldots, K\} \setminus \{ k \} \,\right\}$. The superdifferential of a concave function is defined analogously. Both are referred to as generalized differentials in what follows. The length of the interval corresponding to $\partial_{y_k} f(\bm{y})$ is denoted by 
\begin{equation}
  \lvert \partial_{y_k} f(\bm{y}) \rvert \coloneqq \sup_{a,b \in \partial_{y_k} f(\bm{y})} \; \lvert a - b \rvert.
\end{equation}
If a function $f_{y_k}$ exists such that $f_{y_k}(\bm{y}) \in \partial_{y_k} f(\bm{y}) \; \forall \bm{y} \in \mathcal{Y}$, then $f_{y_k}$ is called a partial generalized derivative of $f$ with respect to $y_k$. The set of all partial generalized derivatives, $f_{y_k}$, is denoted by
$\partial_{y_k} f$.

\subsection{Underlying stochastic process}
\label{ssec:underlying_stochastic_process}

Let $\bigl( X_n \bigr)_{n \geq 1}$ be a discrete-time sto\-chas\-tic process with values in $(\Omega_X, \Filter_X)$. The joint distribution of $\bigl( X_n \bigr)_{n \geq 1}$ on the cylinder set
\begin{equation}
  (\Omega_X^{\mathbb{N}}, \Filter_X^{\mathbb{N}}) \coloneqq \left(  \prod_{n \geq 1} \Omega_X, \prod_{n \geq 1} \Filter_X  \right)
\end{equation}
is denoted by $\Pbb$, the conditional or marginal distributions of an individual random variable $X$ on $(\Omega_X, \Filter_X)$ by $P$ and the natural filtration \cite[Definition~2.32]{Capasso_Bakstein_2015} of the process $\bigl( X_n \bigr)_{n \geq 1}$ by $\bigl( \Filter_X^n \bigr)_{n \geq 1}$. In order to balance generality and tractability, the analysis in this paper is limited to stochastic processes that satisfy the following three assumptions.
\begin{description}
  \item[Assumption 1:] The process $\bigl( X_n \bigr)_{n \geq 1}$ admits a time-homogeneous Markovian representation. That is, there exists a $(\Omega_\Theta, \Filter_\Theta)$-valued stochastic process $\bigl( \Theta_n \bigr)_{n \geq 0}$ adapted to $\bigl( \Filter_X^n \bigr)_{n \geq 1}$ such that
  \begin{align}
    \Pbb(X_{n+1} \in \mathcal{E} \mid X_1 = x_1, \ldots, X_n = x_n) &= \Pbb(X_{n+1} \in \mathcal{E} \mid \Theta_n = \theta_n) \notag \\ 
    &\eqqcolon P_{\theta_n}(\mathcal{E}). 
    \label{def:markov_property}
  \end{align}
  for all $n \geq 1$ and all $\mathcal{E} \in \Filter_X$. The distribution of $X_1$ is denoted by $P_{\theta_0}$, where $\theta_0$ is assumed to be deterministic and known \emph{a priori}. An extension to randomly initialized $\theta_0$ should not be hard but will not be entered.
   
  \item[Assumption 2:] There exists a function $\xi\colon \Omega_\Theta \times \Omega_X \rightarrow \Omega_\Theta$ that is measurable with respect to $\Filter_\Theta \otimes \Filter_X$ and that satisfies 
  \begin{equation}
    \Theta_{n+1} = \xi(\Theta_n,X_{n+1}) \eqqcolon \xi_{\Theta_n}(X_{n+1})
    \label{eq:update_process_statistic}
  \end{equation}
  for all $n \geq 0$. 

  \item[Assumption 3:] For all $\theta \in \Omega_\Theta$, the probability measure $P_{\theta}$, defined in Assumption~1, admits a density $p_{\theta}$ with respect to some $\sigma$-finite reference measure $\mu$.    
\end{description}
The set of distributions $\Pbb$ on $(\Omega_X^{\mathbb{N}}, \Filter_X^{\mathbb{N}})$ that satisfy these three assumptions is denoted by $\mathbb{M}$. The set of distributions $P$ on $(\Omega_X, \Filter_X)$ that admit densities with respect to $\mu$ is denoted by $\mathcal{M}_\mu$.

The above assumptions are rather mild and are introduced primarily to simplify the presentation of the results. In general, the sufficient statistic $\Theta$ can be chosen as a sliding window of past samples, that is, $\Theta_n = (X_{n-m}, \ldots, X_n)$, where $m$ is a finite positive integer. Hence, the presented results apply to every discrete-time Markov process of finite order. However, in order to implement the test in practice, $\Omega_\Theta$ should be sufficiently low-dimensional (compare the examples in Section~\ref{sec:examples}). As long as the existence of the corresponding densities is guaranteed, the reference measure $\mu$ in Assumption~3 can be chosen arbitrarily. This aspect can be exploited to simplify the numerical design of minimax sequential tests and is discussed in more detail in Sections~\ref{sec:discussion} and \ref{sec:examples}.

\subsection{Uncertainty model and hypotheses}
\label{ssec:uncertainty_model_and_hypotheses}

For general Markov processes the question of how to model distributional uncertainty is non-trivial and has far-reaching implications on the definition of minimax robustness. In the most general case the joint distribution $\Pbb$ is subject to uncertainty. However, defining meaningful uncertainty models for $\Pbb$ is an intricate task and is usually neither feasible nor desirable. An approach that is more tractable and more useful in practice is to assume that at any given time instant, $n \geq 1$, the \emph{marginal} or \emph{conditional} distribution of $X_n$ is subject to uncertainty. 

In this paper it is assumed that the conditional distributions $P_\theta$, as defined in \eqref{def:markov_property}, are subject to uncertainty. More precisely, for each $\theta \in \Omega_\Theta$ the conditional distribution $P_{\theta}$ is replaced by an \emph{uncertainty set} of feasible distributions $\Pcal_{\theta} \subset \mathcal{M}_\mu$. This model induces an uncertainty set for $\Pbb$ which is given by 
\begin{equation}
  \Pcal \coloneqq \left\{ \Pbb \in \mathbb{M} : \Pbb = \prod_{n \geq 0} P_{\theta_n}, \; P_{\theta_n} \in \Pcal_{\theta_n} \right\}
  \label{eq:uncertainty_set}
\end{equation}
and is completely specified by the corresponding family of uncertainty sets for the conditional distributions $\{ \Pcal_{\theta} : \theta \in \Omega_\Theta\}$. 

The goal of this paper is to characterize and design minimax optimal sequential tests for multiple hypotheses under the assumption that under each hypothesis the distribution is subject to the type of uncertainty introduced above. That is, each hypothesis is given by
\begin{equation}
  \mathcal{H}_k \colon \Pbb \in \Pcal_k, \quad k = 1, \ldots, K,
  \label{eq:hypotheses}
\end{equation}
where all $\Pcal_k$ are of the form \eqref{eq:uncertainty_set} and are defined by a corresponding family of conditional uncertainty sets $\{ \Pcal_{\theta}^{(k)} : \theta \in \Omega_\Theta\}$. Note that the parameter $\theta$, which corresponds to the sufficient statistic in Assumption~2, does not depend on $k$, that is, the statistic needs to be chosen such that it is sufficient under all hypotheses. Finally, the set $\Pcal_0$ defines the uncertainty in the distribution under which the expected run length is supposed to be minimum. It is assumed to be of the form \eqref{eq:uncertainty_set} as well. In practice, one is often interested in minimizing the worst case expected run length under one hypothesis or under any hypotheses, that is, $\Pcal_0 = \Pcal_k$ for some $k \in \{1, \ldots, K\}$ or $\Pcal_0 = \bigcup_{k=1}^K \Pcal_k$, respectively. In principle, however, $\Pcal_0$ can be chosen freely by the test designer.

Before proceeding, it is useful to illustrate the assumptions on the underlying stochastic process and the proposed uncertainty model with an example. Consider an exponentially weighted moving-average process, that is,
\begin{equation}
  X_{n+1} = \sum_{l=1}^{\infty} a^l X_{n+1-l} + W_{n+1},
\end{equation}
where $a \in (-1,1)$ is a known scalar and $\bigl( W_n \bigr)_{n \geq 1}$ is a sequence of independent random variables that are identically distributed according to $P_W$. This process can equivalently be written as
\begin{equation}
  X_{n+1} = a \Theta_n + W_{n+1},
\end{equation}
where the sufficient statistic $\Theta_n$ can be updated recursively via
\begin{equation}
  \Theta_{n+1} = \xi_{\Theta_n}(X_{n+1}) = \Theta_n + X_{n+1}. 
\end{equation}
In order to introduce uncertainty, it is assumed that with probability $\varepsilon$ the increment $W_n$ is replaced by an arbitrarily distributed outlier. This model yields the following family of conditional uncertainty sets
\begin{equation}
  \Pcal_{\theta} = \left\{ P \in \mathcal{M}_\mu : P(\mathcal{E}) = (1-\varepsilon) P_W(\mathcal{E}-a\theta) + \varepsilon H(\mathcal{E}), \; H \in \mathcal{M}_\mu \right\},
\end{equation}
where $\mathcal{E} \in \Filter_X$ and $\mathcal{E}-a\theta$ is shorthand for $\{ x \in \Omega_X : x+a\theta \in \mathcal{E} \}$. In Section~\ref{sec:examples}, a variant of this example is used to illustrate the design of a minimax optimal test with dependencies in the underlying stochastic process.

\subsection{Testing policies and test statistics}
\label{ssec:Stopping_Rules_decision_rules_test_statistics}

A sequential test is specified via two sequences of randomized decision rules, $\bigl( \psi_n \bigr)_{n \geq 1}$ and $\bigl( \bm{\delta}_n \bigr)_{n \geq 1}$, that are adapted to the filtration $\bigl( \Filter_X^n \bigr)_{n \geq 1}$. Each $\psi_n \colon \Omega_X^n \to [0,1]$ denotes the probability of stopping at time instant $n$. Each $\bm{\delta}_n \colon \Omega_X^n \to \Delta^K$ is a $K$-dimensional vector, $\bm{\delta}_n = (\delta_{1,n}, \ldots, \delta_{K,n})$, whose $k$th element denotes the probability of deciding for $\mathcal{H}_k$, given that the test has stopped at time instant $n$. The randomization is assumed to be performed by independently drawing from a Bernoulli distribution with success probability $\psi_n$ and a discrete distribution on $\{1, \ldots, K\}$ with associated probabilities $(\delta_{1,n}, \ldots, \delta_{K,n})$, respectively. The set of randomized $K$-dimensional decision rules defined on $(\Omega_X^n,\Filter_X^n)$ is denoted by $\Delta_n^{K}$. The stopping time of the test is denoted by $\tau = \tau(\psi)$. 

For the sake of a more concise notation, let $\pi = \bigl( \pi_n \bigr)_{n \geq 1}$ with $\pi_n = (\psi_n, \bm{\delta}_n) \in \Delta_n^1 \times \Delta_n^K$ denote a sequence of tuples of stopping and decision rules. In what follows, $\pi$ is referred to as a \emph{testing policy}, and the set of all feasible policies is denoted by $\Pi \coloneqq \bigtimes_{n \geq 1} \left( \Delta_n^1 \times \Delta_n^K \right)$. 

A test statistic is a stochastic process, $\bigl( T_n \bigr)_{n \geq 0}$, that is adapted to the filtration $\bigl( \Filter_X^n \bigr)_{n \geq 1}$ and allows the stopping and decision rules to be defined as functions mapping from the codomain of $T_n$ to the unit interval. Of particular importance for this paper is the case where the sequence of test statistics, $\bigl( T_n \bigr)_{n \geq 0}$, is itself a time-homogeneous Markov process and the stopping and decision rules are independent of the time index $n$. The corresponding testing policies are in the following referred to as time-homogeneous. This property is formalized in the following definition.
\begin{definition}
  A policy $\pi \in \Pi$ is referred to as \emph{time-homogeneous} if there exists a $(\Omega_T, \Filter_T)$-valued stochastic process $\bigl( T_n \bigr)_{n \geq 0}$ that is adapted to the filtration $\bigl( \Filter_X^n \bigr)_{n \geq 1}$ and it holds that
  \begin{equation}\label{eq:time-homogeneous_policy}
    \psi_n = \psi(T_n) \qquad \text{and} \qquad \bm{\delta}_n = \bm{\delta}(T_n),
  \end{equation}
  where the functions $\psi \colon \Omega_T \to [0,1]$ and $\bm{\delta} \colon \Omega_T \to [0,1]^K$ are independent of the index $n$. 
  \label{def:time-homogeneous_policy}
\end{definition}
Focusing on time-homogeneous policies significantly simplifies both the derivation and the presentation of the main results. Moreover, it will become clear in the course of the paper that optimal tests for time-homogeneous Markov processes can always be realized with time-homogeneous policies so that this restriction can be made without sacrificing optimality.

\subsection{Performance metrics and problem formulation}
\label{ssec:performance_metrics_and_problem_formulation}

The performance metrics considered in this paper are the probability of erroneously rejecting the $k$th hypothesis, $\alpha_k$, and the expected run length of the sequential test, $\gamma$. Both are defined as functions of the testing policy and the true distribution:
\begin{align}
  \gamma(\pi,\Pbb)   &\coloneqq E_{\pi,\Pbb}\bigl[\, \tau(\psi) \,\bigr], 
  \label{eq:runlength} \\
  \alpha_k(\pi,\Pbb) &\coloneqq E_{\pi,\Pbb}\bigl[\, 1-\delta_{k,\tau} \,\bigr], 
  \label{eq:errors}
\end{align}
with $k = 1, \ldots, K$. Here, $E_{\pi,\Pbb}$ denotes the expected value taken jointly with respect to the distribution of $\bigl( X_n \bigr)_{n \geq 1}$ and the randomized policy $\pi$. A generalization to performance metrics that are defined in terms of the pairwise error probabilities is possible but would considerably complicate notation while adding little conceptual insight.

It is important to note that for the design of robust sequential tests the error probabilities and the expected run length need to be treated as equally important performance metrics. On the one hand, reducing the sample size is typically the reason for using sequential tests in the first place. On the other hand, a test whose error probabilities remain bounded over a given uncertainty set, but whose expected run length can increase arbitrarily, cannot be considered robust. In other words, a robust test should not be allowed to delay a decision indefinitely in order to avoid making a wrong decision. 

The first optimality criterion considered in this paper is the weighted sum cost, that is,
\begin{equation}
  L_{\bm{\lambda}}(\pi,\pmb{\Pbb}) = \gamma(\pi,\Pbb_0) + \sum_{k=1}^K \lambda_k \alpha_k(\pi,\Pbb_k),
  \label{eq:weighted_sum_cost}
\end{equation}
where $\pmb{\Pbb} = (\Pbb_0, \ldots, \Pbb_K)$ denotes a $K+1$ dimensional vector of distributions and $\bm{\lambda} = (\lambda_1, \ldots, \lambda_K)$  denotes a $K$ dimensional vector of non-negative cost coefficients. The minimax problem corresponding to the cost function in 
\eqref{eq:weighted_sum_cost} reads as
\begin{equation}
  \adjustlimits \inf_{\pi \in \Pi} \sup_{\pmb{\Pbb} \in \bm{\Pcal}} \; L_{\bm{\lambda}}(\pi,\pmb{\Pbb}),
  \label{eq:unconstrained_minimax}
\end{equation}
where $\pmb{\Pbb} \in \bm{\Pcal}$ is used as a compact notation for $\Pbb_k \in \Pcal_k$,  $k = 0, \ldots, K$. 

The second optimality criterion is the expected run length under constraints on the error probabilities. The corresponding minimax problem reads as
\begin{equation}
  \adjustlimits\inf_{\pi \in \Pi} \sup_{\Pbb_0 \in \Pcal_0} \gamma(\pi,\Pbb_0) 
  \quad \text{s.t.} \quad 
  \sup_{\Pbb_k \in \Pcal_k} \alpha_k(\pi,\Pbb_k) \leq \overline{\alpha}_k,
  \label{eq:constrained_minimax}
\end{equation}
where the constraint holds for all $k = 1, \ldots, K$ and $\overline{\alpha}_k$ denotes an upper bound on the probability of erroneously deciding against $\mathcal{H}_k$. The notation for the minimax optimal policies is fixed below and concludes the section. 

\begin{definition}
  The set of time-homogeneous policies that are optimal in the sense of \eqref{eq:unconstrained_minimax} is denoted by $\PiLambda^*(\bm{\Pcal})$. The set of time-homogeneous policies that are optimal in the sense of \eqref{eq:constrained_minimax} is denoted by $\PiAlpha^*(\bm{\Pcal})$.
  \label{def:minimax_time-homogeneous_policies}
\end{definition}

\section{Optimal Tests}
\label{sec:optimal_tests}

Assume that the distributions $\Pbb_0, \ldots, \Pbb_K \in \mathbb{M}$ are given and fixed. In this case, problems \eqref{eq:unconstrained_minimax} and \eqref{eq:constrained_minimax} reduce to the design of an optimal test for $K$ simple hypotheses, that is,
\begin{equation}
  \inf_{\pi \in \Pi} \; L_{\bm{\lambda}}(\pi,\pmb{\Pbb})
  \label{eq:unconstrained_optimal_test}
\end{equation}
and
\begin{equation}
  \inf_{\pi \in \Pi} \gamma(\pi,\Pbb_0) \quad \text{s.t.} \quad \alpha_k(\pi,\Pbb_k) \leq \overline{\alpha}_k.
  \label{eq:constrained_optimal_test}
\end{equation}
The notation for the corresponding optimal policies is fixed in the next definition.
\begin{definition}
  The set of time-homogeneous policies that are optimal in the sense of \eqref{eq:unconstrained_optimal_test} is denoted by $\PiLambda^*(\pmb{\Pbb})$. The set of time-homogeneous policies that are optimal in the sense of \eqref{eq:constrained_optimal_test} is denoted by $\PiAlpha^*(\pmb{\Pbb})$.
  \label{def:optimal_time-homogeneous_policies}
\end{definition}

The solutions to both the unconstrained problem \eqref{eq:unconstrained_optimal_test} and the constrained problem \eqref{eq:constrained_optimal_test} can be found in the literature. The binary case ($K = 2$) was treated in \cite{Fauss2015} under the same assumptions as stated in Section~\ref{ssec:underlying_stochastic_process}. In \cite{Novikov2009_multiple_hypotheses}, the general solution for an arbitrary number of hypotheses and arbitrary underlying stochastic processes is derived. 

For easier reference the solution of \eqref{eq:unconstrained_optimal_test} is restated in this section. To this end, the functions $\Gfunc, \gfunc \colon \Rnng^{K+1} \to \Rnng$ are introduced. Let
\begin{equation}
  \Gfunc(\z) \coloneqq \sum_{i=1, i \neq k}^K \lambda_i z_i = \left( \sum_{i=1}^K \lambda_i z_i \right) - \lambda_k z_k
  \label{eq:G_function}
\end{equation}
and
\begin{equation}
  \gfunc(\z) \coloneqq \min_{k = 1, \ldots, K} \Gfunc(\z) = \left( \sum_{i=1}^K \lambda_i z_i \right) - \max_{k = 1, \ldots, K} \lambda_k z_k, 
  \label{eq:g_function}
\end{equation}
where $\z = (z_0, \ldots, z_K) \in \Rnng^{K+1}$ and $\bm{\lambda} \in \Rnng^K$ is the vector of cost coefficients introduced in Section~\ref{ssec:performance_metrics_and_problem_formulation}. Note that both $\Gfunc$ and $\gfunc$ are independent of $z_0$; defining them as functions of the $K+1$ dimensional vector $\z$ unifies the notation in what follows.

The cost function that characterizes the optimal test is stated in the following theorem. It extends Theorem~2.1 in \cite{Fauss2015} to multiple hypotheses.
\begin{theorem}
  Let $\bm{\lambda} \geq 0$, let $\pmb{\Pbb} \in \mathbb{M}^{K+1}$, and let $\rhofunc \colon \rhodomain \to \Rnng$, where
  \begin{equation}
    \rhodomain \coloneqq \Rnng^{K+1} \times \Omega_\Theta.
  \end{equation}
  The integral equation
  \begin{equation}
    \rhofunc(\z, \theta) = \min \left\{ \gfunc(\z) \,,\, z_0 + \int \rhofunc \bigl( \z \p_{\theta}(x), \xi_\theta(x) \bigr) \, \mu(\dint x) \right\},
    \label{eq:optimal_rho}
  \end{equation}
  with $\gfunc$ defined in \eqref{eq:g_function}, has a unique solution and it holds that
  \begin{equation}
    \inf_{\pi \in \Pi} \; L_{\bm{\lambda}}(\pi,\pmb{\Pbb}) = \rhofunc(\bm{1},\theta_0).
  \end{equation}
  \label{th:optimal_cost_function}
\end{theorem}
Theorem~\ref{th:optimal_cost_function} follows directly from Theorem~5 and Lemma~6 in \cite{Novikov2009_multiple_hypotheses} and the Markov property of the stochastic process $\bigl( X_n \bigr)_{n \geq 1}$. Also, compare Theorem~5 in \cite{Fauss2016_thesis}. The optimal test statistic and testing policies are obtained by comparing the cost for stopping with the expected cost for continuing under the optimal policy. 
\begin{corollary} 
  The optimal test statistic of a test solving \eqref{eq:unconstrained_optimal_test} is given by
  \begin{equation}
    T_n(x_1,\ldots,x_n) = (\z_n, \theta_n) \in \rhodomain, 
    \label{def:optimal_test_statistic}
  \end{equation}
  where $\theta_n$ is a sufficient statistic for $(x_1,\ldots,x_n)$ in the sense of \eqref{def:markov_property} and $\z_n = (z_{0,n}, \ldots, z_{K,n})$ is the vector of likelihood ratios (Radon--Nikodym derivatives)
  \begin{equation}
    z_{k,n} = \prod_{i=1}^n \frac{\mathrm{d} P_{\theta_{i-1}}^{(k)}}{\mathrm{d} \mu}(x_i) = \prod_{i=1}^n p_{\theta_{i-1}}^{(k)}(x_i). 
    \label{def:z}
  \end{equation}
  By Assumption~2 in Section~\ref{ssec:underlying_stochastic_process}, the test statistic can be calculated recursively via
    \begin{align}
    \theta_{n+1} &= \xi_{\theta_n}(x_{n+1}), & \theta_0 &\colon \text{given a priori}  \\[1ex]
    \z_{n+1} &= \z_n \p_{\theta_n}(x_{n+1}), & \z_0 &= \bm{1},
  \end{align}
  with $\xi_\theta$ is defined in \eqref{eq:update_process_statistic}.
  \label{col:optimal_test_statistic}
\end{corollary}
\begin{corollary}
  Let $\bm{\lambda} \geq 0$, let $\pmb{\Pbb} \in \mathbb{M}^{K+1}$, and let $\rhofunc$ be as defined in Theorem~\ref{th:optimal_cost_function}. A policy $\pi$ is time-homogeneous in the sense of Definition~\ref{def:time-homogeneous_policy} and optimal in the sense of \eqref{eq:unconstrained_optimal_test}, that is, $\pi \in \PiLambda^*(\pmb{\Pbb})$, if and only if its stopping and decision rules are of the form
  \begin{equation}
    \psi_n = \psi(\z_n, \theta_n) \qquad \text{and} \qquad \delta_{k,n} = \delta_k(\z_n),
    \label{eq:optimal_time-homogenoues_policies}
  \end{equation}
  where
  \begin{gather}
    \mathcal{I}\bigl(\{\gfunc(\z) < \rhofunc(\z, \theta)\}\bigr) \leq \psi(\z, \theta) \leq \mathcal{I}\bigl(\{\gfunc(\z) \leq \rhofunc(\z, \theta)\}\bigr), \\[1ex]
    \delta_k(\z) \leq \mathcal{I}\left(\left\{\Gfunc(\z) = \gfunc(\z)\right\}\right),
  \end{gather}
  and $(\z_n, \theta_n)$ is defined in Corollary~\ref{col:optimal_test_statistic}.
  \label{col:optimal_policies}
\end{corollary}

Corollary~\ref{col:optimal_test_statistic} and Corollary~\ref{col:optimal_policies} follow immediately from Theorem~6 in \cite{Novikov2009_multiple_hypotheses} and characterize the set of optimal time-homogeneous policies under Markov assumptions. The two components of the optimal test statistic $(\z_n, \theta_n)$ correspond to the two types of information that are necessary to apply the optimal stopping rule. The likelihood ratios $\z$ are needed to evaluate the cost for stopping, while the state of the Markov process $\theta$ is needed to evaluate the conditional expectation that determines the cost for continuing. 

The policies defined in Corollary~\ref{col:optimal_policies} are optimal in the sense of the unconstrained problem \eqref{eq:unconstrained_optimal_test}. The solution of the constrained problem is closely related, but its statement is deferred to the next section since it relies on properties of the optimal cost function $\rhofunc$ that need to be established first. Before turning to the latter, a more compact notation and a useful way of characterizing the performance of time-homogeneous tests of the form \eqref{def:optimal_test_statistic} are introduced.

For policies of the form given in Corollaries~\ref{col:optimal_test_statistic} and \ref{col:optimal_policies}, it is convenient to define the expected run length and the error probabilities of the underlying test as functions of the initial state of the test statistic, that is,
\begin{align}
  \gamma_{\pi, \Pbb}(\z, \theta) &\coloneqq E_{\pi,\Pbb} \bigl[\, \tau(\psi) \mid \bm{Z}_0 = \z, \Theta_0 = \theta \,\bigr]
  \label{eq:conditional_runlength} \\
  \alpha_{\pi, \Pbb}^{(k)}(\z, \theta) &\coloneqq E_{\pi,\Pbb}\bigl[\, 1-\delta_{\tau,k} \mid \bm{Z}_0 = \z, \Theta_0 = \theta \,\bigr],  
  \label{eq:conditional_errors}
\end{align}
where $\gamma_{\pi,\Pbb} \colon \rhodomain \to \Rnng \cup \{ \infty \}$ and $\alpha_{\pi,\Pbb} \colon \rhodomain \to [0,1]$. Since $\bigl(\bm{Z}_n, \Theta_n\bigr)_{n \geq 0}$ is a time-homogeneous Markov process, $\gamma_{\pi,\Pbb}$ and $\alpha_{\pi,\Pbb}^{(k)}$ solve the Chapman--Kolmogorov equations \cite{EoM1990}
\begin{equation}
  \gamma_{\pi,\Pbb}(\z,\theta) = (1-\psi(\z,\theta)) \left( 1 + \int \gamma_{\pi,\Pbb} \bigl( \z \p_\theta(x), \xi_\theta(x) \bigr) \, P_\theta(\dint x) \right)
  \label{eq:ck_runlength}
\end{equation}
and
\begin{equation}
  \alpha_{\pi,\Pbb}^{(k)}(\z,\theta) = \psi(\z,\theta) \, \bigl(1-\delta_k(\z,\theta)\bigr) + \bigl( 1-\psi(\z,\theta) \bigr) \int \alpha_{\pi,\Pbb}^{(k)} \bigl( \z \p_\theta(x), \xi_\theta(x) \bigr) \, P_\theta(\dint x). 
  \label{eq:ck_errors}
\end{equation}

Conditioning on the true initial states of the statistics reduces the conditional performance metrics in \eqref{eq:conditional_runlength} and \eqref{eq:conditional_errors} to the unconditional ones in \eqref{eq:runlength} and \eqref{eq:errors}, that is,
\begin{equation}
  \gamma_{\pi,\Pbb}(\bm{1},\theta_0) = \gamma(\pi,\Pbb)
  \qquad \text{and} \qquad 
  \alpha_{\pi,\Pbb}^{(k)}(\bm{1},\theta_0) = \alpha_k(\pi,\Pbb).
  \label{eq:conditional_unconditional_performace}
\end{equation}
Being able to characterize the performance of a test in terms of solutions of \eqref{eq:ck_runlength} and \eqref{eq:ck_errors} is of critical importance for the proofs given in later sections.

In order to simplify the notation of the central integral equations of this section, let $\bigl\{ \mu_{\z,\theta} : (\z,\theta) \in \rhodomain \bigr\}$ and $\bigl\{ P_{\z,\theta} : (\z,\theta) \in \rhodomain \bigr\}$ be two families of probability measures on $\Omega_\rho$ that are defined via
\begin{align} 
  \mu_{\z,\theta}(\mathcal{E}_{\z} \times \mathcal{E}_{\theta}) &\coloneqq \mu \bigl( \left\{\, x \in \Omega_X : \z \p_{\theta}(x) \in \mathcal{E}_{\z} \,,\, \xi_{\theta}(x) \in \mathcal{E}_{\theta} \right\} \,\bigr),
  \label{eq:mu_z_theta} \\[1ex]
  P_{\z,\theta}(\mathcal{E}_{\z} \times \mathcal{E}_{\theta}) &\coloneqq P_{\theta} \bigl( \left\{\, x \in \Omega_X : \z \p_{\theta}(x) \in \mathcal{E}_{\z} \,,\, \xi_{\theta}(x) \in \mathcal{E}_{\theta} \right\} \,\bigr), 
  \label{eq:P_z_theta}
\end{align}
where $\mathcal{E}_z \times \mathcal{E}_{\theta} \in \rhofilter$, with $\rhofilter$ denoting the natural $\sigma$-algebra on $\rhodomain$. The notation $P_{\z,\theta}^{(k)}$ is used to refer to the probability measure in \eqref{eq:P_z_theta} with $P_{\theta}$ chosen to be $P_{\theta}^{(k)}$. Using this notation, the integral equations in \eqref{eq:optimal_rho}, \eqref{eq:ck_runlength} and \eqref{eq:ck_errors} can be written more compactly as
\begin{equation}
  \rhofunc = \min \left\{ \gfunc \,,\, z_0 + \int \rhofunc \, \dint \mu_{\z, \theta} \right\}
\end{equation}
and
\begin{gather}
  \gamma_{\pi,\Pbb} = (1-\psi) \left( 1 + \int \gamma_{\pi,\Pbb} \, \dint P_{\z,\theta} \right),
  \label{eq:ck_runlength_compact} \\
  \alpha_{\pi,\Pbb}^{(k)} = \psi \, \bigl(1-\delta_k\bigr) + \bigl( 1-\psi \bigr) \int \alpha_{\pi,\Pbb}^{(k)} \, \dint P_{\z,\theta}.
  \label{eq:ck_errors_compact}
\end{gather}
Both notations are used in what follows.

\section{\texorpdfstring{Properties of the Cost Function $\rho_{\lambda}$}{Properties of the Cost Function}}
\label{sec:properties_of_rho}

While Theorem~\ref{th:optimal_cost_function} is well-known in the literature, the cost function $\rhofunc$ in \eqref{eq:optimal_rho} has rarely been studied in detail. The connection between $\rhofunc$ and the properties of sequential tests using policies $\pi \in \PiLambda^*(\pmb{\Pbb})$ is the subject of this section. It extends and generalizes the results in Section 3 of \cite{Fauss2015}.

\begin{theorem} 
  Let $\rhofunc$ be as defined in Theorem~\ref{th:optimal_cost_function}, and let $\PiLambda^*(\pmb{\Pbb})$ be as defined in Definition~\ref{def:optimal_time-homogeneous_policies}. For all $\pi \in \PiLambda^*(\pmb{\Pbb})$, it holds that
  \begin{equation}
    \rhofunc(\z, \theta) = z_0 \, \gamma_{\pi,\Pbb_0}(\z, \theta) + \sum_{k=1}^K \lambda_k z_k \, \alpha_{\pi,\Pbb_k}^{(k)}(\z, \theta), 
    \label{eq:rho_sum}
  \end{equation}
  where $\gamma_{\pi,\Pbb}$ and $\alpha_{\pi,\Pbb}^{(k)}$ are defined in \eqref{eq:conditional_runlength} and \eqref{eq:conditional_errors}, respectively.
  \label{th:rho_sum}
\end{theorem}

See Appendix~\ref{apd:proof_rho_sum} in \cite{Fauss2018_supplement} for a proof. Theorem~\ref{th:rho_sum} connects the performance metrics in Section~\ref{ssec:performance_metrics_and_problem_formulation} to the optimal cost function $\rhofunc$ and is key to solving the unconstrained minimax problem \eqref{eq:unconstrained_minimax}. However, obtaining a solution to the constrained minimax problem \eqref{eq:constrained_minimax} additionally requires control over the individual error probabilities. In what follows, a connection between the latter and the partial generalized derivatives of $\rhofunc$ is established. First, some useful technical properties of $\rhofunc$ are shown.
\begin{lemma}
  For all $\bm{\lambda} \geq 0$ and all $\pmb{\Pbb} \in \mathbb{M}^{K+1}$, the function $\rhofunc$ that solves \eqref{eq:optimal_rho} is non-decreasing, concave and homogeneous of degree one in $\z$.
  \label{lm:rho_properties}
\end{lemma}
A proof is detailed in Appendix~\ref{apd:proof_rho_properties} in \cite{Fauss2018_supplement}. Lemma~\ref{lm:rho_properties} is significant for two reasons. First, it ensures that $\rhofunc$ admits a generalized differential, a property that is used in the next theorem to establish a connection between $\rhofunc$ and the error probabilities $\alpha_k$. Second, being concave and homogeneous qualifies $\rhofunc$ as an $f$-dissimilarity, that is, a statistical measure for the joint similarity of $P_{\theta}^{(0)}, \ldots, P_{\theta}^{(K)}$. This property of $\rhofunc$ is not used explicitly in the derivations of the minimax optimal sequential test, but it puts the results into context and will later be shown to provide a unified interpretation of minimax optimal sequential and minimax optimal fixed sample size test. A more detailed discussion of this aspect is deferred to Section~\ref{sec:discussion}.

\begin{theorem}
  Let $\pmb{\Pbb} \in \mathbb{M}^{K+1}$ and let $\PiLambda^*(\pmb{\Pbb})$ be as defined in Definition~\ref{def:optimal_time-homogeneous_policies}. For $\rhofunc$ as defined in Theorem~\ref{th:optimal_cost_function}, and $\gamma_{\pi,\Pbb}, \alpha_{\pi,\Pbb}^{(k)}$ as defined in \eqref{eq:conditional_runlength} and \eqref{eq:conditional_errors}, respectively, it holds that: 
  \begin{enumerate}
    \item for all $\pi \in \PiLambda^*(\pmb{\Pbb})$ and all $k = 1, \ldots, K$
          \begin{align}
            \gamma_{\pi,\Pbb_0} &\in \partial_{z_0} \rhofunc, \\
            \lambda_k \alpha_{\pi,\Pbb_k}^{(k)} &\in \partial_{z_k} \rhofunc.
          \end{align}
    \item for all $\pi \in \PiLambda^*(\pmb{\Pbb})$, all $(\z,\theta) \in \rhodomain$, and all $k = 1, \ldots, K$
          \begin{align}
            \Big\{ \gamma_{\pi,\Pbb_0}(\z, \theta) : \pi \in \PiLambda^*(\pmb{\Pbb}) \Big\} &= \partial_{z_0} \rhofunc(\z, \theta), 
            \label{eq:policy_diff_0} \\
            \Big\{ \lambda_k \alpha_{\pi,\Pbb_k}^{(k)}(\z, \theta) : \pi \in \PiLambda^*(\pmb{\Pbb}) \Big\} &= \partial_{z_k} \rhofunc(\z, \theta). 
            \label{eq:policy_diff_k}
          \end{align}
  \end{enumerate}
  \label{th:rho_derivatives}
\end{theorem}
Theorem~\ref{th:rho_derivatives} is proven in Appendix~\ref{apd:proof_rho_derivatives} in \cite{Fauss2018_supplement}. Its two parts correspond to a \emph{global} and a \emph{local} statement about the generalized differentials of $\rhofunc$. The first part states that for all $\pi \in \PiLambda^*(\pmb{\Pbb})$, the functions $\gamma_{\pi,\Pbb_0}$ and $\alpha_{\pi,\Pbb_k}^{(k)}$ are valid generalized differentials of $\rhofunc$. The second part states that at every point $(\z,\theta) \in \rhodomain$, the generalized differential of $\rhofunc$ coincides with the set of all error probabilities that can be realized by policies $\pi \in \PiLambda^*(\pmb{\Pbb})$. Note that since the optimal policy is deterministic on the interior of the decision regions, the remaining degrees of freedom in terms of the error probabilities are exclusively due to the randomization on the boundary of the stopping and decision regions. That is, the subdifferentials on the right-hand sides of \eqref{eq:policy_diff_0} and \eqref{eq:policy_diff_k} are spanned by policies with different randomizations on the left-hand side. If only deterministic policies are allowed, the result only holds for the end points of the interval, that is, the left and right derivative.

The local statement in Theorem~\ref{th:rho_derivatives} cannot be extended to a global statement since the integral equations \eqref{eq:set_integral_diff_0} and \eqref{eq:set_integral_diff_k} in \cite{Fauss2018_supplement} establish a coupling between the local differentials so that they cannot be chosen independently of each other. This coupling reflects the fact that changing the randomization on the boundaries of a decision region also affects the overall performance of the corresponding sequential test. 

Based on Theorem~\ref{th:rho_derivatives}, the following optimality result can be obtained for the constrained sequential testing problem in \eqref{eq:constrained_optimal_test}.
\begin{theorem}
  Let $\pmb{\Pbb} \in \mathbb{M}^{K+1}$, and let $\pi^* \in \PiAlpha^*(\pmb{\Pbb})$, that is, $\pi^*$ solves \eqref{eq:constrained_optimal_test}. If
  \begin{equation}
    \bm{\lambda}^* \in \argmax_{\bm{\lambda} \in \Rnng^K} \; \left\{ \rhofunc(\bm{1}, \theta_0) - \sum_{k=1}^K \lambda_k \overline{\alpha}_k \right\},
    \label{eq:lambda_opt}
  \end{equation}
  with $\rhofunc$ defined in Theorem~\ref{th:optimal_cost_function}, then for all $\pi \in \Pi_{\lambda^*}^*(\pmb{\Pbb})$ and all $k = 1, \ldots, K$ it holds that
  \begin{align}
    \bigl\lvert \gamma(\pi,\Pbb_0) - \gamma(\pi^*,\Pbb_0) \bigr\rvert &\leq \bigl\lvert \partial_{z_0} \rho_{\bm{\lambda}^*}(\bm{1}, \theta_0) \bigr\rvert, 
    \label{eq:constrained_optimal_gamma} \\[1ex]
    \lambda_k^* \, \bigl\lvert \alpha_k(\pi,\Pbb_k) - \overline{\alpha}_k \bigr\rvert &\leq \bigl\lvert \partial_{z_k} \rho_{\bm{\lambda}^*}(\bm{1}, \theta_0) \bigr\rvert,
    \label{eq:constrained_optimal_alpha}
  \end{align}
  Moreover, it holds that
  \begin{equation}
    \Pi_{\lambda^*}^*(\pmb{\Pbb}) \cap \PiAlpha^*(\pmb{\Pbb}) \neq \emptyset,
  \end{equation}
  that is, there exists at least one $\pi \in \Pi_{\lambda^*}^*(\pmb{\Pbb})$ that is optimal in the sense of \eqref{eq:constrained_optimal_test}.
  \label{th:optimal_constrained_test}
\end{theorem}
A proof is detailed in Appendix~\ref{apd:proof_optimal_constrained_test} in \cite{Fauss2018_supplement}. The significance of Theorem~\ref{th:optimal_constrained_test} lies in the fact that it characterizes solutions of the constrained problem \eqref{eq:constrained_optimal_test} in terms of the optimal cost function of the unconstrained problem \eqref{eq:unconstrained_optimal_test}. From an algorithmic point of view, Theorem~\ref{th:optimal_constrained_test} makes it possible to design constrained sequential tests via a systematic optimization of the cost coefficients $\bm{\lambda}$ instead of Monte Carlo simulations or resampling techniques \cite{Tartakovsky2003, Sochman2005}. 

Theorems~\ref{th:rho_sum}, \ref{th:rho_derivatives} and \ref{th:optimal_constrained_test} form the basis for the derivation of minimax optimal tests, which will be characterized as optimal tests for least favorable distributions. The latter are introduced and discussed in the next section.

\section{Least Favorable Distributions}
\label{sec:least_favorable_distributions}

The counterpart of the optimal testing problem investigated in the previous sections is the problem of determining the least favorable distributions for a given testing policy $\pi$. In this case the unconstrained problem \eqref{eq:unconstrained_optimal_test} reduces to
\begin{equation}
  \sup_{\pmb{\Pbb} \in \bm{\Pcal}} \; L_{\bm{\lambda}}(\pi, \pmb{\Pbb}) = \sup_{\pmb{\Pbb} \in \bm{\Pcal}} \left( \gamma(\pi,\Pbb_0) + \sum_{k=1}^K \lambda_k \alpha_k(\pi,\Pbb_k) \right).
  \label{eq:lfds_joint}
\end{equation}
Since the expected run length and the error probabilities are coupled only via the policy, the joint problem in \eqref{eq:lfds_joint} decouples into $K+1$ individual maximization problems
\begin{equation}
  \sup_{\Pbb_0 \in \Pcal_0} \gamma(\pi,\Pbb_0) 
  \qquad \text{and} \qquad
  \sup_{\Pbb_k \in \Pcal_k} \alpha_k(\pi, \Pbb_k), \quad k = 1, \ldots, K,
  \label{eq:lfd_decoupled}
\end{equation}
which can be solved independently.

For arbitrary stopping and decision rules, solving the problems in \eqref{eq:lfd_decoupled} exactly is challenging and, in general, the least favorable distributions depend on the time instant $n$ as well as on the history of the random process. However, for time-homogeneous policies of the form \eqref{eq:optimal_time-homogenoues_policies}, a more elegant solution can be obtained.
\begin{theorem}
  Let $\Pcal$ be an uncertainty set of the form \eqref{eq:uncertainty_set}, and let $\gamma_{\pi,\Pcal} \colon \rhodomain \to \Rnng \cup \{ \infty \}$ and $\alpha_{\pi,\Pcal}^{(k)} \colon \rhodomain \to \Rnng$. For all testing policies $\pi$ of the form \eqref{eq:optimal_time-homogenoues_policies} and all $k = 1, \ldots, K$, it holds that the integral equations
  \begin{gather}
    \gamma_{\pi,\Pcal} = (1-\psi) \left( 1 + \sup_{H \in \Pcal_{\theta}} \int \gamma_{\pi,\Pcal}\bigl( \z \p_\theta(x), \xi_\theta(x) \bigr) \, H(\dint x) \right) 
    \label{eq:lfds_equation_runlength} \\
    \alpha_{\pi,\Pcal}^{(k)} = \psi(1-\delta_k) + (1-\psi) \left( \sup_{H \in \Pcal_{\theta}} \int \alpha_{\pi,\Pcal}^{(k)} \bigl( \z \p_\theta(x), \xi_\theta(x) \bigr) \, H(\dint x) \right) 
    \label{eq:lfds_equation_errors}
  \end{gather}
  have unique solutions.
  \label{th:lfds_equations}
\end{theorem}
\begin{theorem}
   Let $\bm{\Pcal} = (\Pcal_0, \ldots, \Pcal_K)$ be uncertainty sets of the form \eqref{eq:uncertainty_set}, and let $\pi$ be of the form \eqref{eq:optimal_time-homogenoues_policies}.
   \begin{itemize}
      \item If it holds that for every $(\z,\theta) \in \rhodomain$
            \begin{equation}
              \Qcal_{\z,\theta}^{(0)} \coloneqq \argmax_{H \in \Pcal_{\theta}^{(0)}} \int \gamma_{\pi, \mathcal{P}_0}\bigl( \z \p_\theta(x), \xi_\theta(x) \bigr) \, H(\dint x) \neq \emptyset
            \end{equation}
            with $\gamma_{\pi, \mathcal{P}}$ defined in \eqref{eq:lfds_equation_runlength}, then every distribution 
            \begin{equation}
            \Qbb_0 \in \Qcal_0 \coloneqq \left\{ \Pbb \in \mathbb{M} : \Pbb = \prod_{n \geq 0} P_{\z_n,\theta_n} \,,\, P_{\z_n,\theta_n} \in \Qcal^{(0)}_{\z_n,\theta_n} \right\}
            \end{equation}
            is least favorable with respect to the expected run length of the test, that is,
            \begin{equation}
              \gamma(\pi,\Qbb_0) = \sup_{\Pbb \in \Pcal_0} \gamma(\pi,\Pbb).
            \end{equation}

     \item If it holds that for every $(\z,\theta) \in \rhodomain$ 
           \begin{equation}
             \Qcal_{\z,\theta}^{(k)} \coloneqq \argmax_{H \in \Pcal_{\theta}^{(k)}} \int \alpha_{\pi, \mathcal{P}_k}^{(k)} \bigl( \z \p_\theta(x), \xi_\theta(x) \bigr) \, H(\dint x) \neq \emptyset
           \end{equation}
           with $k = 1, \ldots, K$ and $\alpha_{\pi, \mathcal{P}}^{(k)}$ defined in \eqref{eq:lfds_equation_errors}, then every distribution 
           \begin{equation}
             \Qbb_k \in \Qcal_k \coloneqq \left\{ \Pbb \in \mathbb{M} : \Pbb = \prod_{n \geq 0} P_{\z_n,\theta_n} \,,\, P_{\z_n,\theta_n} \in \Qcal^{(k)}_{\z_n,\theta_n} \right\}
           \end{equation}
           is least favorable with respect to the probability of erroneously deciding against $\mathcal{H}_k$, that is,
           \begin{equation}
             \alpha_k(\pi,\Qbb_k) = \sup_{\Pbb \in \Pcal_k} \alpha_k(\pi,\Pbb).
           \end{equation}
  \end{itemize}  
  \label{th:lfds}
\end{theorem}
Theorems~\ref{th:lfds_equations} and \ref{th:lfds} are proven in Appendices~\ref{apd:proof_lfds_equations} and \ref{apd:proof_lfds} in \cite{Fauss2018_supplement}, respectively. Note that the distributions whose densities are used to update the test statistic $\z \p_\theta$ are not the least favorable distribution, but are determined by the poilcy $\pi$; compare Corollary~\ref{col:optimal_test_statistic}. From Theorem~\ref{th:lfds} it follows that, under the least favorable distributions, the process $\bigl(X_n\bigr)_{n \geq 1}$ is a Markov process with sufficient statistic $(\z_n,\theta_n)$. That is, the least favorable distributions adapt to the policy of the test as well as to the history of the stochastic process. 

It is worth highlighting that even in the case where the $X_n$ are independent, that is, no additional statistic $\theta_n$ is required, the least favorable distributions still generate a Markov process with sufficient statistic $\z$. That is, the least favorable distributions depend on the history of the process via the likelihood ratios. This aspect distinguishes robust sequential hypothesis tests from existing robust tests for fixed sample sizes or for change points. It is discussed in more detail in Section~\ref{sec:discussion}.

Having characterized optimal tests and least favorable distributions, everything is in place for the derivation of minimax optimal sequential tests.

\section{Minimax Optimal Sequential Tests}
\label{sec:minimax_optimal_sequential_tests}

In this section sufficient conditions for strict minimax optimality of sequential hypothesis tests are given. Following the procedure for the optimal sequential test without uncertainty, the solutions of the unconstrained problem \eqref{eq:unconstrained_minimax} are derived first and are then shown to contain a solution of the constrained problem \eqref{eq:constrained_minimax}. 

The following three theorems are stated in sequence and constitute the main contribution of the paper. A discussion and an interpretation of the results is deferred to the next section.  

\begin{theorem}
  Let $\bm{\lambda} \geq 0$, let $\bm{\Pcal} = (\Pcal_0, \ldots, \Pcal_K)$ be uncertainty sets of the form \eqref{eq:uncertainty_set}, and let $\rhofunc, \dfunc \colon \rhodomain \to \Rnng$ and $\Dfunc \colon \rhodomain \times \mathcal{M}_{\mu} \to \Rnng$. The equation system
  \begin{align}
    \rhofunc(\z, \theta) &= \min \left\{\, \gfunc(\z) \,,\, z_0 + \dfunc(\z, \theta) \,\right\}   
    \label{eq:minimax_rho} \\[1.5ex]
    \dfunc(\z, \theta) &= \sup_{\bm{P} \in \bm{\Pcal}_{\theta}} \; \Dfunc(\z, \theta; \bm{P}) 
    \label{eq:minimax_d} \\[1ex]
    \Dfunc(\z, \theta; \bm{P}) &= \int \rhofunc \bigl( \z \p_\theta(x), \xi_\theta(x) \bigr) \, \mu(\dint x)
    \label{eq:minimax_D}
  \end{align}
  with $\gfunc$ defined in \eqref{eq:g_function} has a unique solution.
  \label{th:minimax_equations}
\end{theorem}

\begin{theorem}
  Let $\bm{\lambda} \geq 0$, let $\bm{\Pcal} = (\Pcal_0, \ldots, \Pcal_K)$ be uncertainty sets of the form \eqref{eq:uncertainty_set}, and let $\rhofunc, \dfunc$, and $\Dfunc$ be as defined in Theorem~\ref{th:minimax_equations}. If for all $(\z,\theta) \in \rhodomain$
  \begin{equation}
    \bm{\Qcal}_{\z,\theta} \coloneqq \argmax_{\bm{P} \in \bm{\Pcal}_{\theta}} \; \Dfunc(\z, \theta; \bm{P}) \neq \emptyset,
    \label{eq:sup=max}
  \end{equation}
  then every policy $\pi \in \PiLambda^*(\pmb{\Qbb})$ with 
  \begin{equation}
    \pmb{\Qbb} \in \pmb{\Qcal} = \left\{ \pmb{\Pbb} \in \mathbb{M}^{K+1} : \Pbb_k = \prod_{n \geq 0} P_{\z_n,\theta_n}^{(k)} \,,\, \bm{P}_{\z_n,\theta_n} \in \bm{\Qcal}_{\z_n,\theta_n} \right\}
  \end{equation}
  is minimax optimal in the sense of \eqref{eq:unconstrained_minimax}, that is, 
  \begin{equation}
    \bigl\{\, \pi \in \PiLambda^*(\pmb{\Qbb}) : \pmb{\Qbb} \in \pmb{\Qcal} \,\bigr\} \subset \PiLambda^*(\pmb{\Pcal }).
  \end{equation}
  \label{th:minimax_tests}
\end{theorem}

\begin{theorem}
  Let $\pmb{\Pcal} = (\Pcal_0, \ldots, \Pcal_K)$ be uncertainty sets of the form \eqref{eq:uncertainty_set}, and let $\pi^* \in \PiAlpha^*(\pmb{\Pcal})$ such that $(\pi^*, \pmb{\Pbb}^*)$ solves \eqref{eq:constrained_minimax}. If
  \begin{equation}
    \bm{\lambda}^* \in \argmax_{\bm{\lambda} \in \Rnng^K} \; \left\{ \rhofunc(\bm{1}, \theta_0) - \sum_{k=1}^K \lambda_k \overline{\alpha}_k \right\}
    \label{eq:lambda_minimax}
  \end{equation}
  with $\rhofunc$ defined in Theorem~\ref{th:minimax_equations}, then for all $\pi \in \Pi_{\lambda^*}^*(\pmb{\Qbb})$ with $\pmb{\Qbb} \in \pmb{\Qcal}$ defined in Theorem~\ref{th:minimax_tests} and all $k = 1, \ldots, K$ it holds that
  \begin{align}
    \bigl\lvert \gamma(\pi,\Qbb_0) - \gamma(\pi^*,\Pbb_0^*) \bigr\rvert &\leq \bigl\lvert \partial_{z_0} \rho_{\bm{\lambda}^*}(\bm{1}, \theta_0) \bigr\rvert, 
    \label{eq:constrained_minimax_gamma} \\[1ex]
    \lambda_k^* \, \bigl\lvert \alpha_k(\pi,\Qbb_k) - \overline{\alpha}_k \bigr\rvert &\leq \bigl\lvert \partial_{z_k} \rho_{\bm{\lambda}^*}(\bm{1}, \theta_0) \bigr\rvert.    \label{eq:constrained_minimx_alpha}
  \end{align}
  Moreover, it holds that
  \begin{equation}
    \bigl\{\, \pi \in \Pi_{\bm{\lambda}^*}^*(\pmb{\Qbb}) : \pmb{\Qbb} \in \pmb{\Qcal} \,\bigr\} \cap \PiLambda^*(\pmb{\Pcal }) \neq \emptyset,
  \end{equation}
  that is, there exists at least one pair $(\pi, \pmb{\Qbb})$, with $\pmb{\Qbb} \in \pmb{\Qcal}$ and $\pi \in \Pi_{\bm{\lambda}^*}^*(\pmb{\Qbb})$, that is optimal in the sense of \eqref{eq:constrained_minimax}.
  \label{th:minimax_constrained_test}
\end{theorem}
Theorems~\ref{th:minimax_equations}, \ref{th:minimax_tests}, and \ref{th:minimax_constrained_test} are proven in Appendices~\ref{apd:proof_minimax_equations}, \ref{apd:proof_minimax_tests}, and \ref{apd:proof_constrained_minimax_test} in \cite{Fauss2018_supplement}, respectively.

\section{Discussion}
\label{sec:discussion}

Theorems~\ref{th:minimax_tests} and \ref{th:minimax_constrained_test} in the previous section provide sufficient conditions for the characterization of minimax optimal tests in terms of optimal testing policies and least favorable distributions. In this section the question of existence is discussed, and an interpretation in terms of statistical similarity measures is given that provides additional insight and establishes a connection to minimax optimal fixed samples size tests and change detection procedures.

\subsection{Statistical Similarity Measures}
\label{ssec:statistical_similarity_measures}

In order to obtain a better conceptual understanding of minimax optimal sequential tests, it is helpful to introduce a class of statistical similarity measures known as $f$-dissimilarities. They were first proposed by Gy\"orfi and Nemetz \cite{GyorfiNemetz1975, GyorfiNemetz1977, GyorfiNemetz1978} as an extension of $f$-divergences to multiple distributions and play an important role in the theory of statistical decision making. In particular, the connection between $f$-dissimilarities and Bayesian risks has been a topic of high interest in statistics \cite{Nguyen2009}, signal processing \cite{Varshney2011} and machine learning \cite{Reid2011}. 

In this section it is shown that the function $\rhofunc$ in Theorem~\ref{th:minimax_equations} induces an $f$-dissimilarity and that this $f$-dissimilarity provides a sufficient characterization of the minimax optimal test. For this purpose a variation on the concept of $f$-dissimilarities is useful, which is defined as follows:
\begin{definition}
  Let $P_1,\ldots,P_K$ be probability measures on a measurable space $(\Omega,\Filter)$, and let $f\colon \Rnng^K \times \Omega \to \mathbb{R}$, with $f = f(\bm{y}, \omega) = f(y_1, \ldots, y_K, \omega)$, be homogeneous of degree one and concave in $(y_1,\ldots,y_K)$. The functional
  \begin{equation}
    I_f(P_1,\ldots,P_K) = \int f(p_1(\omega), \ldots, p_K(\omega), \omega) \, \mu(\mathrm{d} \omega)
    \label{eq:f-simmilarity}
  \end{equation}
  is called $f$-similarity of $P_1,\ldots,P_K$.
  \label{def:f_similarity}
\end{definition}
$I_f$ in \eqref{eq:f-simmilarity} is referred to as a \emph{similarity} measure since for concave and homogeneous functions $f$, the functinal $-I_f = I_{-f}$ is a \emph{dissimilarity} measure in the sense of Gy\"orfy and Nemetz \cite[Definition 1]{GyorfiNemetz1978}. Allowing $f$ to depend on the integration variable directly is a minor generalization that, in the context of this paper, allows the similarity measure to depend on the history of the random process. Similar generalizations have been introduced in the literature before \cite{Rockafellar1968, Papageorgiou1997, Breuer2016}.

Using Definition~\ref{def:f_similarity}, an intuitive characterization of the family of least favorable distributions can be given in terms of a corresponding family of $f$-similarities.
\begin{corollary}
  At every time instant $n \geq 1$, the least favorable distributions of $X_{n+1}$, conditioned on the state $(\bm{Z}_n,\Theta_n) = (\z,\theta)$, are the feasible distributions that are most similar with respect to the $f$-similarity defined by
  \begin{equation}
    f_{\z,\theta}(\bm{y}, x) = \rhofunc\bigl( \z \bm{y} \,,\, \xi_\theta(x) \bigr),
  \end{equation}
  with $\rhofunc$ given in Theorem~\ref{th:minimax_equations}.
  \label{cor:lfds_similarity}
\end{corollary}
The family of $f$-divergences, defined by $\{ f_{\z,\theta} : (\z,\theta) \in \rhodomain \}$, can be interpreted as follows: The equation system in Theorem~\ref{th:minimax_equations} defines the optimal cost function $\rhofunc$ which in turn defines the similarity measure in Corollary~\ref{cor:lfds_similarity}. The sequential aspect of the test is captured by the parameters $(\z,\theta)$. The likelihood ratios $\z$ determine the weights of the individual densities. That is, the larger $z_k$, the larger the influence of the $k$th distribution on the similarity measure. In terms of the underlying hypothesis test, a high value of $z_k$ implies that the test is likely to decide in favor of $\mathcal{H}_k$. Consequently, depending on whether $\mathcal{H}_k$ is true or not, the least favorable distributions need to be as similar or dissimilar to $\mathcal{H}_k$ as possible. On the other hand, a low value of $z_k$ implies that a decision in favor of $\mathcal{H}_k$ is highly unlikely so that the corresponding distribution contributes little to the overall similarity measure. The influence of the parameter $\theta$, that is, the history of the underlying stochastic process, does not affect the relative weighting of the individual distributions, but rather the shape of $\rhofunc$. It strongly depends on the stochastic process itself, and it is less predictable than the effect of $\z$. 

The use of statistical similarity measures for the design of (sequential) tests for composite hypotheses has been suggested previously in the literature. However, as discussed in the introduction, in most works a suitable similarity measure is chosen \emph{beforehand}---usually based on asymptotic results, bounds, or approximations---and is used as a \emph{surrogate} objective whose optimization is easier than solving the actual testing problem. Here, by contrast, it is shown that the minimax testing problem \emph{induces} a similarity measure and that the distributions that maximize the latter solve the former.

It is also instructive to compare the minimax optimal tests in the previous section to tests whose least favorable distributions achieve Choquet capacities. In \cite{Huber1973, Fauss2016_old_bands}, it is shown that in this case the pair of least favorable distributions jointly minimizes \emph{all} $f$-divergences over the uncertainty sets. This implies that the least favorable distributions do not depend on the decision rule and can be calculated \emph{a priori}; compare also \cite{Verdu1984}. For the minimax sequential test, this is no longer the case. Instead, the least favorable distributions need to minimize a particular $f$-dissimilarity that depends on $\rhofunc$ and the current state of the test statistic. Interestingly, the existence of capacity achieving distributions seems to be closely related to the existence of a single threshold test. Both fixed sample size detection and quickset change detection are problems for which optimal single threshold procedures exist and whose least favorable distributions have been shown to be of the capacity type \cite{Huber1973, Unnikrishnan2011}. Sequential detection, on the other hand, is inherently a multi-threshold procedure. This difference is subtle but turns out to be crucial. An intuitive explanation is as follows: For a single threshold test the least favorable distributions are those that maximize or minimize the drift of the test statistic towards the threshold. In the example of quickest change detection, maximizing the expected detection delay corresponds to minimizing the drift towards the threshold under the post-change distribution. On the other hand, in order to maximize the expected run length of a two-sided sequential test, one needs to minimize the drift of the test statistic towards \emph{either threshold}. Depending on the data observed so far, this can translate to a drift towards the upper or towards the lower threshold; the least favorable distributions become data dependent. Similar considerations also apply to fixed sample size tests for multiple hypotheses which might explain why, as of now, no capacity based minimax solutions for the latter for can be found in the literature.

In summary, the discussion shows that minimax optimal sequential tests are based the same principles that underpin existing minimax tests. optimality is achieved by using a policy that leads to the \emph{best separation} of the \emph{most similar} distributions. The difference is that in the sequential case the similarity measure depends on the testing policy and on the observed data.

\subsection{Existence}
\label{ssec:existence}

The results presented so far allow for some statements about the existence of minimax optimal tests. Stronger statements can be made for the unconstrained problem formulation \eqref{eq:unconstrained_minimax} than for the constrained formulation \eqref{eq:constrained_minimax}. The former is considered first.

Since by Theorem~\ref{th:optimal_cost_function} it holds that
\begin{align*}
  \adjustlimits \inf_{\pi \in \Pi} \sup_{\pmb{\Pbb} \in \bm{\Pcal}} \; L_{\bm{\lambda}}(\pi,\pmb{\Pbb}) \leq \gfunc(\bm{1}) \leq \sum_{k=1}^K \lambda_k,
\end{align*}
the minimax optimal objective value in \eqref{eq:unconstrained_minimax} is guaranteed to be finite for all cost coefficients $\bm{\lambda} \geq 0$ and all uncertainty sets $\pmb{\Pcal}$. This includes scenarios where two uncertainty sets overlap or are identical. However, it does not imply that a testing policy exists that achieves this value. In order to guarantee \eqref{eq:sup=max}, that is, that the supremum is attained, the uncertainty sets $\Pcal_1, \ldots, \Pcal_K$ need to be compact \cite[Theorem 4.16]{Rudin1976}. While this assumption seems restrictive in theory, it is automatically fulfilled when the problem is discretized in order to be solved numerically. Moreover, given that the least favorable distributions are maximizers of concave functionals, we conjecture that compactness is not necessary for the existence of minimax solutions. A more thorough discussion would need to go into the technical details of optimization on function spaces, which is beyond the scope of this paper.

For the constrained problem \eqref{eq:constrained_minimax} the situation is more involved, since, in order to be able to satisfy the constraints on the error probabilities, the uncertainty sets $\Pcal_1, \ldots, \Pcal_K$ need to be sufficiently \emph{separated}, i.e., the distance between the two sets needs to be large enough to statistically separate them using a finite number of samples. The appropriate way to measure this distance is via the $f$-similarity that is induced by the corresponding cost function $\rhofunc$. This leads to the following result.
\begin{corollary}
  The minimax optimal objective in \eqref{eq:constrained_minimax} is finite, if and only if the right hand side of \eqref{eq:lambda_minimax} is bounded, i.e., if
  \begin{equation}
    \sup_{\bm{\lambda} \in \Rnng^K} \; \left\{ \rhofunc(\bm{1}, \theta_0) - \sum_{k=1}^K \lambda_k \overline{\alpha}_k \right\} < \infty,
    \label{eq:lagrange_dual}
  \end{equation}
  with $\rhofunc$ defined in Theorem~\ref{th:minimax_equations}.
  \label{col:constrained_minimax_existence}
\end{corollary}
A proof of Corollary~\ref{col:constrained_minimax_existence} is detailed in Appenidx~\ref{apd:proof_col:constrained_minimax_existence}. As in the unconstrained case, existence of a finite supremum does not imply that an optimal policy exists, unless the uncertainty sets are compact. Finally, note that no deterministic policy might exist that attains the bounds on the error probabilities with equality, which implies that, in general, the optimal policy is randomized. In view of Theorem~\ref{th:rho_derivatives}, this happens when the bounds are located in the interior of the subdifferentials of the optimal cost function.

Interestingly, the conditions for the existence of minimax optimal sequential tests are rather mild, especially for the unconstrained problem formulation. This is in contrast to the fixed sample size case for which much stricter sufficient conditions are given in the literature \cite{Huber1973, Fauss2016_old_bands}. This suggests the conjecture that these stricter conditions are only necessary to guarantee that the optimal policy and the least favorable distributions are \emph{decoupled} and that minimax optimal fixed sample size tests for arbitrary uncertainty sets exist but that they require a joint design of the policy and the least favorable distributions. Unfortunately, applying the results presented in this paper to the fixed sample size case is not straightforward. First, the concept of ordered samples, which is essential to sequential hypothesis testing and induces state-dependent least favorable distributions, in general does not apply to fixed sample size tests. Second, owing to the deterministic stopping rule, the policies of fixed sample size tests are not time-homogeneous in the sense of Definition~\ref{def:time-homogeneous_policy}. However, the second numerical example given in the next section indicates that it might be possible to implement truncated tests with time-homogeneous policies---see Section~\ref{ssec:binominal_ar1} for more details. 
    
Finally, it should be highlighted that the least favorable distributions do not depend on the reference measure $\mu$. More precisely, if $\mu$ is absolutely continuous with respect to an alternative reference measure $\tilde{\mu}$, it follows from the homogeneity of $\rhofunc$ that 
\begin{align*}
  \Dfunc(\z, \theta; \bm{P}) = \int \rhofunc\bigl(\z \p, \xi_\theta \bigr) \, \dint \mu &= \int \rhofunc\bigl(\z \p, \xi_\theta \bigr) \frac{\dint \mu}{\dint \tilde{\mu}} \, \dint \tilde{\mu} \\
  &= \int \rhofunc\biggl(\z \p \frac{\dint \mu}{\dint \tilde{\mu}}, \xi_\theta \biggr)  \, \dint \tilde{\mu} 
  \eqqcolon \int \rhofunc\bigl(\z \tilde{\p}, \xi_\theta \bigr)  \, \dint \tilde{\mu},
\end{align*}
where $\tilde{\bm{p}} = (\tilde{p}_0, \ldots, \tilde{p}_K)$ are probability densities with respect to $\tilde{\mu}$.

\section{Example and Numerical Results}
\label{sec:examples}

In order to illustrate the presented results, two examples for minimax optimal sequential tests are given in this section. First, a test for three hypotheses that is robust with respect to the error probabilities is designed under the assumption that $\bigl( X_n \bigr)_{n \geq 1}$ is a sequence of independent random variables with identical uncertainty sets. Second, the case when the underlying processes admits dependencies is illustrated by solving the Kiefer--Weiss problem for a variant of the uncertainty model in Section~\ref{ssec:uncertainty_model_and_hypotheses}. In both examples the targeted error probabilities are set to $\overline{\alpha}_k = 0.01$ for all $k = 1, \ldots, K$.

The test design is based on the following iterative procedure. First, all $P_{\theta}^{(k)}$ are initialized to some feasible distribution. Keeping these distributions fixed, an optimal sequential test is designed by solving \eqref{eq:lambda_opt} for $\bm{\lambda}^*$ and $\rho_{\bm{\lambda}^*}$. In the second step, $\rho_{\bm{\lambda}^*}$ is kept fixed and the distributions are updated by solving the optimization problem in \eqref{eq:sup=max}. Both steps are iterated until the changes in the function $\rho_{\bm{\lambda}^*}$ are small enough to assume convergence. For both examples convergence was assumed if the relative difference between two consecutive approximations of $\rho_{\bm{\lambda}^*}$ fell below $10^{-3}$. The question whether this procedure is guaranteed to converge in general is certainly worth investigating but beyond the scope of this paper.

It should be noted that this iterative procedure does \emph{not} alternate between the design of optimal testing policies and least favorable distributions. Unless the procedure has converged, the distributions that solve \eqref{eq:sup=max} are not least favorable in the sense of Theorem~\ref{th:lfds}. Moreover, the test statistic, which is part of the optimal policy, depends on the likelihood ratios and is hence affected by the update of the distributions in the second step. 

In order to solve \eqref{eq:lambda_opt} and \eqref{eq:sup=max} numerically, both the state space $\rhodomain$ and the sample space $\Omega_X$ are discretized using a regularly spaced grid, and linear interpolation is used to evaluate functions between grid points. This straightforward approach works well for the examples presented here. However, if a larger number of hypotheses or more complex dependencies need to be considered, more sophisticated approximations need to be used \cite{Ferrari2005, Kunsch2017}. The linear programming algorithm detailed \cite{Fauss2015} was used to efficiently solve \eqref{eq:lambda_opt} jointly for $\bm{\lambda}^*$ and $\rho_{\bm{\lambda}^*}$. However, in principle, any suitable convex optimization algorithm can be used to solve \eqref{eq:lambda_opt}, and any method for solving nonlinear integral equations can be used to obtain $\rho_{\bm{\lambda}^*}$. 

In both examples, the distributional uncertainty is of the density band type \cite{Kassam1981}, that is,
\begin{equation}
  \Pcal = \bigl\{ P \in \mathcal{M}_{\mu} : p' \leq p \leq p'' \bigr\},
  \label{eq:band_model}
\end{equation}
where $0 \leq p' \leq p'' \leq \infty$, $P'(\Omega_X) \leq 1$ and $P''(\Omega_X) \geq 1$. Here, $P'$ and $P''$ denote the measures corresponding to $p'$ and $p''$, respectively. The reason for using this uncertainty model is twofold. First, it contains several popular uncertainty models as special cases, for example, the $\varepsilon$-contamination model \cite{Huber1965}, the bounded distribution function model \cite{Oesterreicher1978, Hafner1993} and the $f$-divergence ball model \cite{Fauss2018}. Second, an efficient iterative algorithm for the minimization of convex functionals of probability distribution under density band constraints exists that makes it possible to obtain accurate numerical solutions with moderate computational efforts. A more detailed discussion of the band model, its properties and how to obtain the least favorable distributions numerically can be found in \cite{Kassam1981} and \cite{Fauss2016_old_bands}.

\subsection{IID process under three hypotheses}
\label{ssec:iid}

For the first example all $X_n$, $n \geq 1$, are assumed to be independent and distributed on the interval $\Omega_X = [-1,1]$. Let $P_n$ denote the distribution of $X_n$. The task is to decide between the following three hypotheses:
\begin{equation}
  \mathcal{H}_1 \colon P_n \in \Pcal_1, \qquad \mathcal{H}_2 \colon P_n \in \Pcal_2, \qquad \mathcal{H}_3 \colon P_n = \mathcal{U}_{[-1,1]},
\end{equation}
for all $n \geq 1$. Here, $\mathcal{U}_{[a,b]}$ denotes the continuous uniform distribution on the interval $[a, b]$ and the uncertainty sets $\Pcal_1$, $\Pcal_2$ are of the form \eqref{eq:band_model} with
\begin{align}
  p_1'(x) &= a e^{-2x} + 0.1, & p_1''(x) &= a e^{-2x} + 0.3, \\
  p_2'(x) &= a e^{2x} + 0.1,  & p_2''(x) &= a e^{2x} + 0.3,
\end{align}
where $a \approx 0.1907$ was chosen such that $P_1'(\Omega_X) = P_2'(\Omega_X) = 0.9$ and $P_1''(\Omega_X) = P_2''(\Omega_X) = 1.1$. The expected runlength was minimized under $\mathcal{H}_3$, i.e., $P_n = \mathcal{U}_{[-1,1]}$ for all $n$. Moreover, in order to keep the domain of the cost function two-dimensional, the reference measure $\mu$ was set to $\mu = \mathcal{U}_{[-1,1]}$ so that $z_0 = z_3 = 1$ and $\rhofunc$ becomes a function of $(z_1, z_2)$ only. 

The scenario in this example can arise, for instance, in monitoring applications, where $\mathcal{H}_3$ corresponds to an ``in control'' state in which the distribution of the data is known almost exactly, while $\mathcal{H}_1$ and $\mathcal{H}_2$ correspond to two different ``out of control'' states, with only partially known distributions. If it needs to be established that the system is ``in control'' before a certain procedure starts, it is reasonable to minimize the expected run length of the test under the ``in control'' distribution while still requiring it to be insensitive to distributional uncertainties in the ``out of control'' states.

\begin{figure}[tb]
  \centering
  \includegraphics[width=0.5\textwidth]{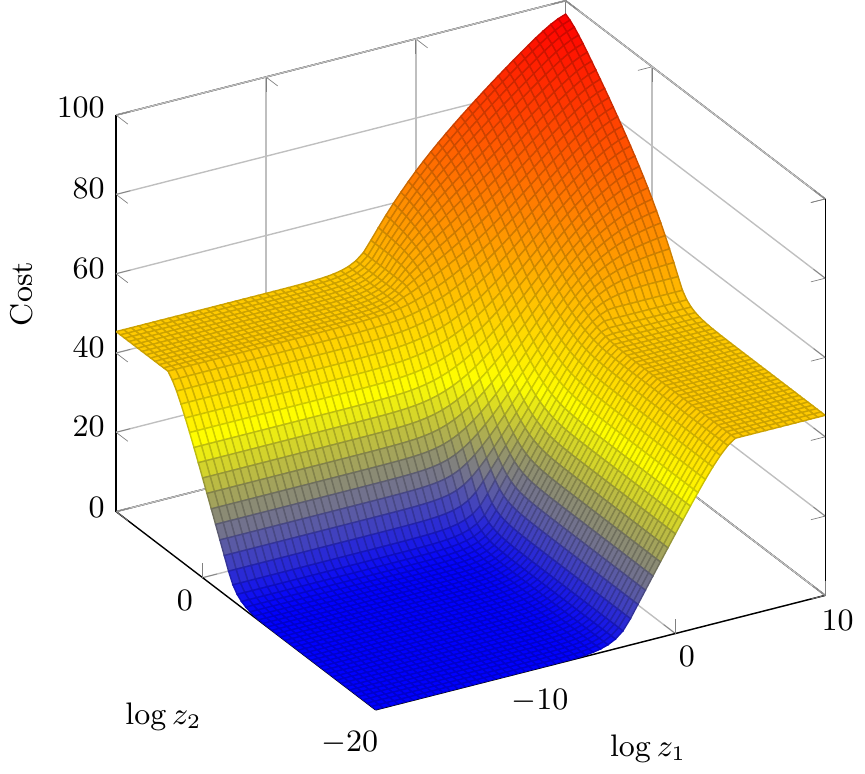}%
  \includegraphics[width=0.5\textwidth]{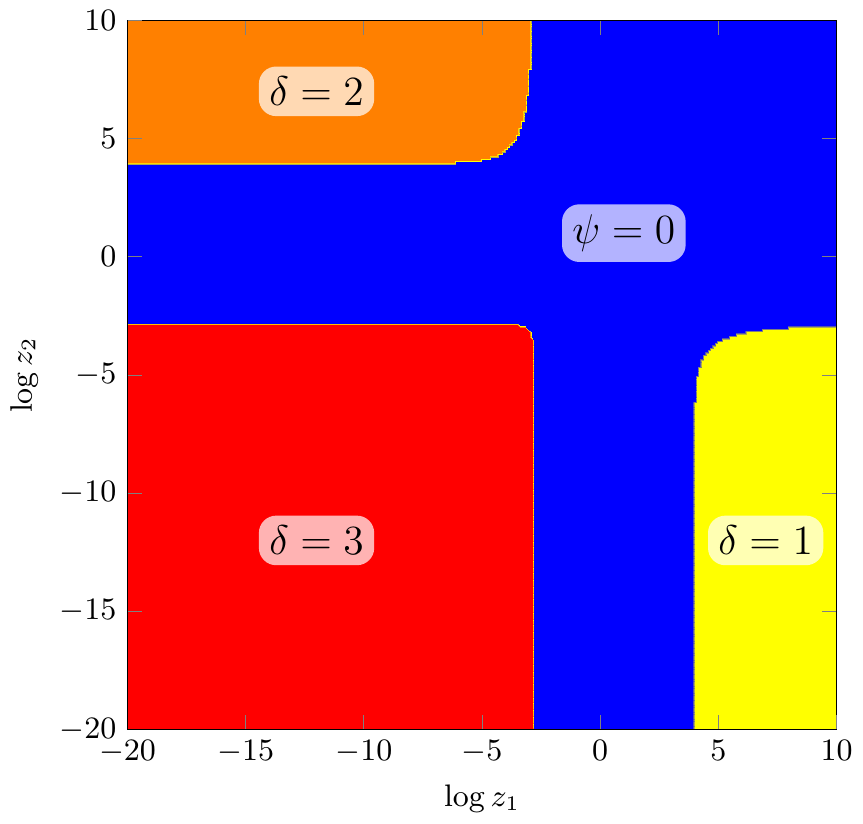}
  \caption{Segment of the optimal cost function $\rhofunc$ in logarithmic scale and the corresponding testing policy as a function of the log-likelihood ratios.}
  \label{fig:opt_policy_1}
\end{figure}

In order to solve this example numerically, the likelihood ratio plane $(z_1, z_2)$ was discretized on $[-20, 10] \times [-20, 10]$ using $301 \times 301$ uniformly spaced grid points, and the sample space $\Omega_X = [-1, 1]$ was discretized using $101$ uniformly spaced grid points. The design procedure detailed above converged after five iterations. The optimal weights were found to be $\bm{\lambda}^* \approx (133.41, 133.41, 45.41)$. The resulting cost function $\rho_{\bm{\lambda}^*}$, as well as the corresponding testing policy, are depicted in Figure~\ref{fig:opt_policy_1}. While the cost function as such provides little insight, the testing policy lends itself to an intuitive interpretation. In analogy to the regular sequential probability ratio test (SPRT), the minimax optimal test consists of two corridors that correspond to a binary test between $\mathcal{H}_{\{1,2\}}$ and $\mathcal{H}_3$, respectively. Interestingly, there is a rather sharp intersection of the two corridors so that the test quickly reduces to a quasi-binary scenario. 

\begin{figure}[tb]
  \centering
  \includegraphics[width=0.5\textwidth]{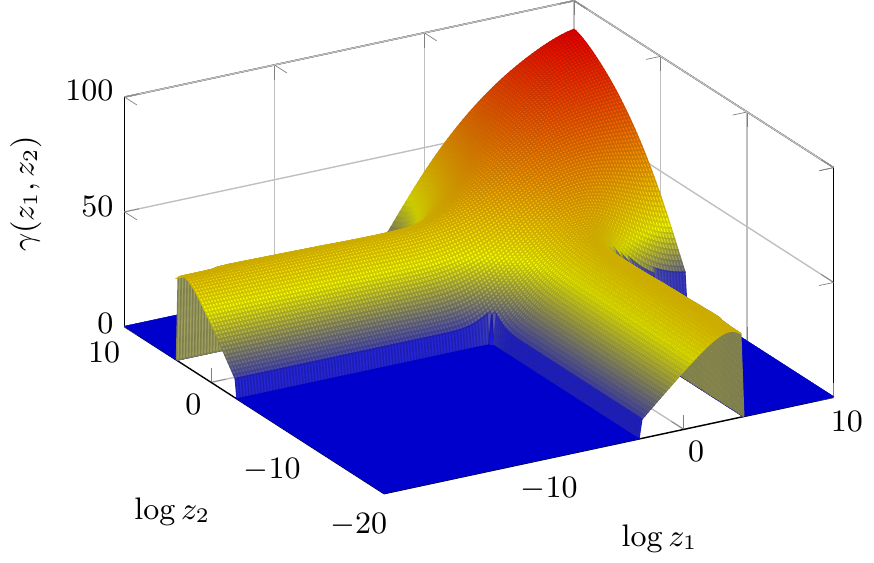}%
  \includegraphics[width=0.5\textwidth]{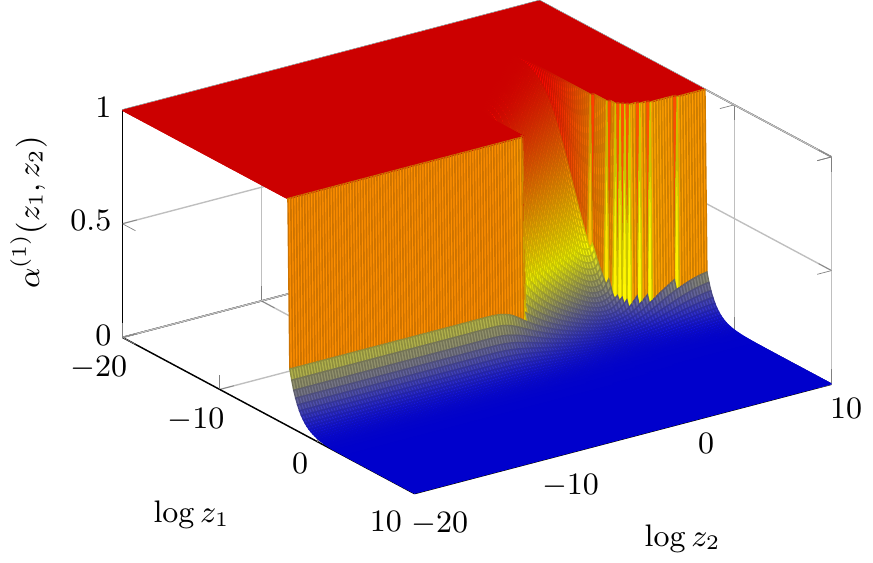}\\
  \includegraphics[width=0.5\textwidth]{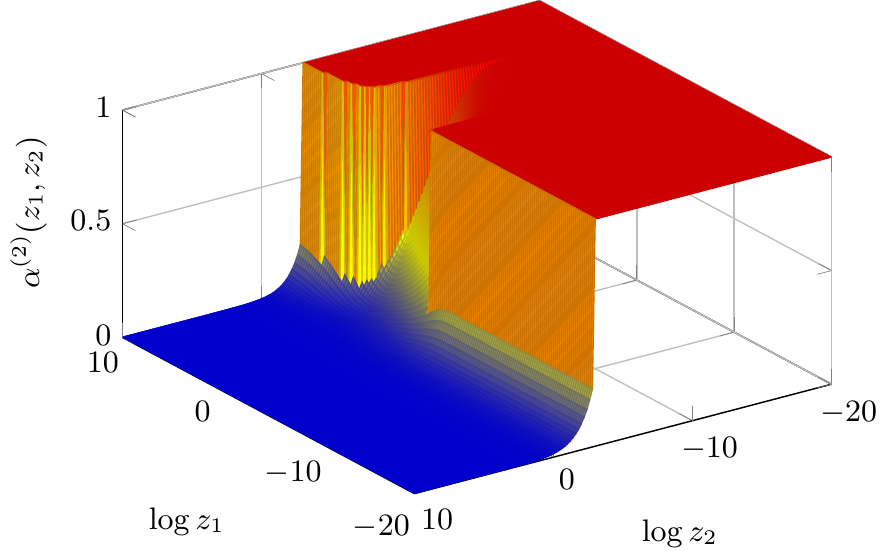}%
  \includegraphics[width=0.5\textwidth]{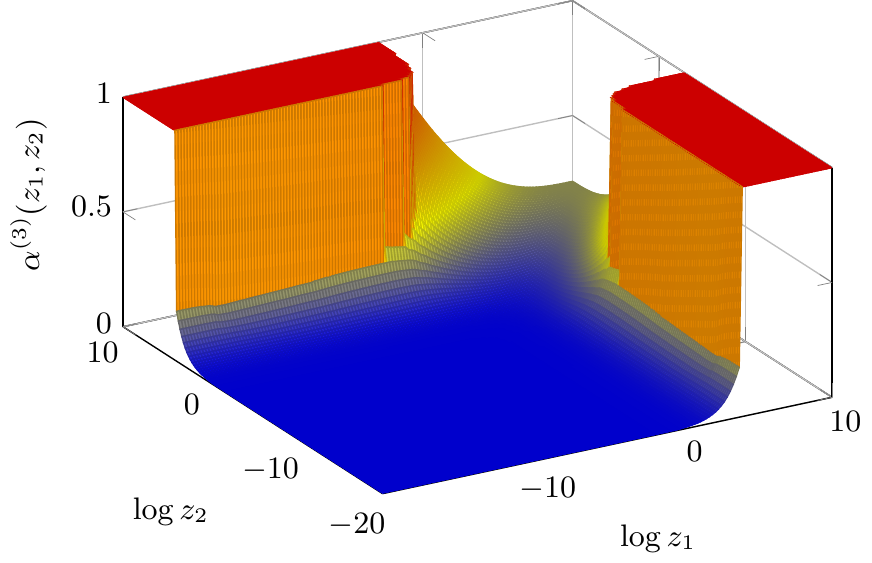}
  \caption{Performance metrics as functions of the log-likelihood ratios. Clockwise from the top left: expected run length, error probability of the first, second and third type.}
  \label{fig:performance_1}
\end{figure}

The expected run length and the error probabilities as functions of the state of the test statistic can be obtained either via the partial derivatives of $\rho_{\bm{\lambda}^*}$ or by solving the integral equations \eqref{eq:conditional_errors} and \eqref{eq:runlength} and are depicted in Figure~\ref{fig:performance_1}. The ``blocky'' appearance of some of the functions is due to them having being downsampled to a coarser grid for plotting. Moreover, no smoothing was applied in order not to smear the hard transitions between the decision regions. Finally, note that the plots are oriented differently to provide a better visual representation of the respective function.

\begin{figure}[tb]
  \centering
  \includegraphics[width=0.5\textwidth]{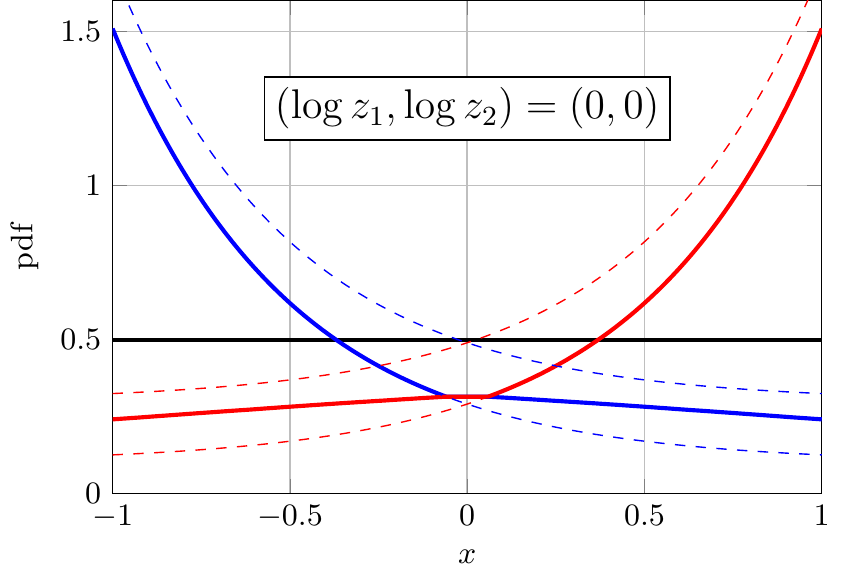}%
  \includegraphics[width=0.5\textwidth]{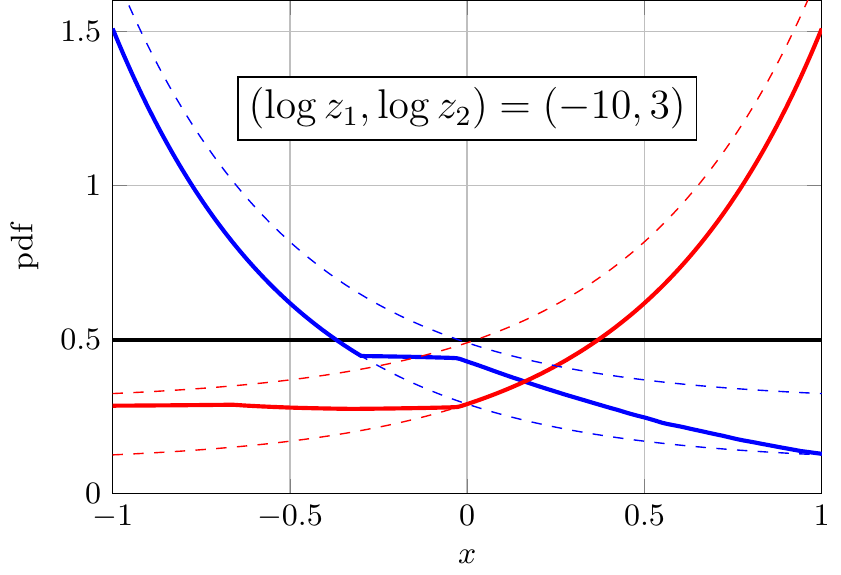}\\
  \includegraphics[width=0.5\textwidth]{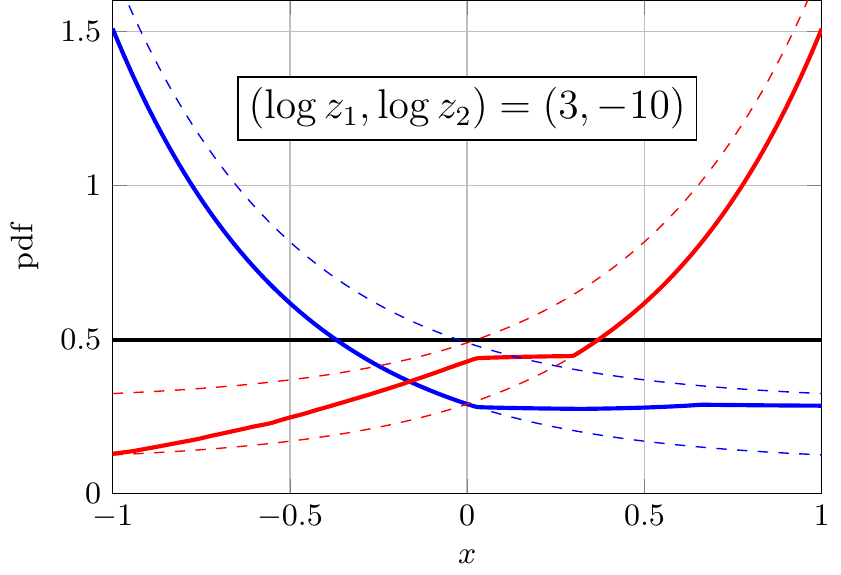}%
  \includegraphics[width=0.5\textwidth]{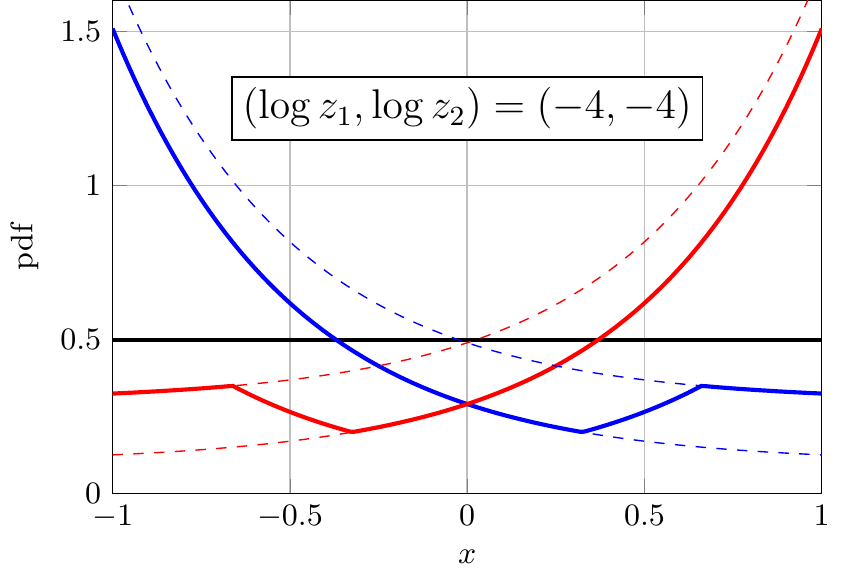}\\
  \includegraphics{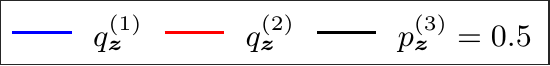}
  \caption{Examples of least favorable distributions for different states of the test statistic.}
  \label{fig:lfds_1}
\end{figure}

The stopping and decision rules in Figure~\ref{fig:opt_policy_1} are depicted as functions of the log-likelihood ratios. The latter are in turn defined in terms of the least favorable distributions, that is, the distributions that solve the maximization in \eqref{eq:sup=max}. Four examples of densities of least favorable distributions are depicted in Figure~\ref{fig:lfds_1}. As can be seen, the densities change significantly, depending on the state of the test statistic. In the top left plot, the test statistic is in its initial state, meaning that there is no preference for either hypothesis. Consequently, the least favorable densities are chosen such that all three distributions are equally similar to each other which in this case implies that they are symmetric around the $y$-axis and that $q_{\z}^{(1)}$ and $q_{\z}^{(2)}$ jointly mimic the uniform density $p_{\z}^{(3)}$. Also note that $q_{\z}^{(1)}$ and $q_{\z}^{(2)}$ overlap on an interval around $x = 0$ so that observations in this interval are statistically indistinguishable under $\mathcal{H}_1$ and $\mathcal{H}_2$. As the test statistic is updated, the least favorable distributions change. In the upper right and the lower left plot of Figure~\ref{fig:lfds_1}, two cases are depicted where the test has a strong preference for $\mathcal{H}_1$ or $\mathcal{H}_2$, respectively; compare the decision regions in Figure~\ref{fig:opt_policy_1}. In both cases the least favorable densities are no longer symmetric, but their probability masses are shifted, their tail-behavior is noticeably different and the interval of overlap can no longer be observed. Finally, in the lower right plot, there is a strong preference for $\mathcal{H}_3$, which leads to $q_{\z}^{(1)}$ and $q_{\z}^{(2)}$ both shifting as much probability mass as possible to their tails in order to reduce the significance of the corresponding observations. It is interesting to observe the effect that an imminent decision for $\mathcal{H}_3$ has on $q_{\z}^{(1)}$ and $q_{\z}^{(2)}$, namely, that they become less similar to each other in order to increase the joint similarity to $p_{\z}^{(3)}$. This is in contrast to the initial state depicted in the upper left plot, where $q_{\z}^{(1)}$ and $q_{\z}^{(2)}$ also try to approximate $p_{\z}^{(3)}$, but at the same time need to be similar to each other as well.

In order to verify the numerical results, \num{10000} Monte Carlo simulations were performed using the testing policy depicted in Figure~\ref{fig:opt_policy_1}. The observations were drawn from the least favorable distributions, which were calculated on the fly by solving \eqref{eq:sup=max} for the current weights $(z_1, z_2)$. The resulting confusion matrix and the average run length of the tests are shown in Table~\ref{table:mc_results_1}.

\begin{table}[tb]
  \centering
  \includegraphics{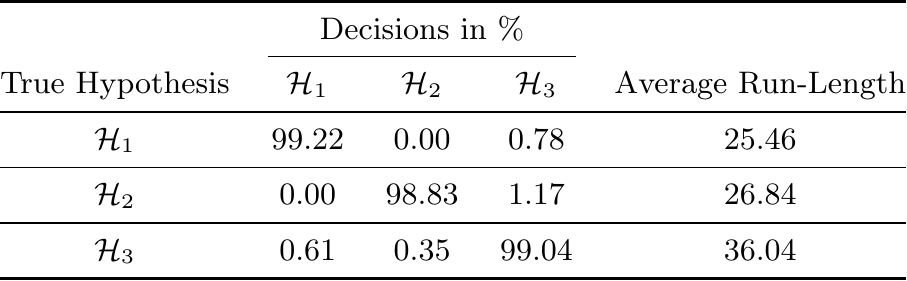}%
  \caption{Results of $10^4$ Monte Carlo runs of the minimax optimal test using the policy depicted in Figure\ref{fig:opt_policy_1}. The target detection probability is $99\%$, the theoretical expected run length under $\mathcal{H}_3$ is $E_{\Pbb_3}\bigl[\, \tau(\psi) \,\bigr] \approx 36.31$ samples.}
  \label{table:mc_results_1}
\end{table}

\subsection{Binomial AR(1) process under two hypotheses}
\label{ssec:binominal_ar1}

A binomial AR(1) process $\bigl( X_n \bigr)_{n \geq 1}$ is a homogeneous Markov process with transition probabilities
\begin{equation}
  P(X_n=m \mid X_{n-1} = \theta) = \smashoperator{\sum_{i = \max\{0 \,,\, m+\theta-M\}}^{\min\{m\,,\,\theta\}}} b_{i,\theta,m,M} \, \beta_0^i (1-\beta_0)^{\theta-i} \beta_1^{m-i} (1-\beta_1)^{M-\theta+i-m},
  \label{eq:binom_ar1}
\end{equation}
where $m = 0, \ldots, M$,
\begin{equation}
  b_{i,\theta,m,M} = \binom{l}{i} \binom{M-\theta}{m-i}
\end{equation}
and $(\beta_0, \beta_1) = \bm{\beta} \in (0,1)^2$ characterizes the dependence structure of the process. See \cite[Definition 1.1, Remark 1.2]{Weiss2013} for a formal definition and more details on the parameter $\bm{\beta}$ and its feasible values. The sufficient statistic of the binomial AR(1) process is given by $\Theta_n = X_n$ with $\Omega_\theta = \Omega_X = \{0, 1, \ldots, M\}$. In what follows, $M = 7$.

Let $P_{\theta,\bm{\beta}}$ denote the distribution in \eqref{eq:binom_ar1}, that is, the conditional distribution of $X_n$ given $X_{n-1} = \theta$. The following two simple hypotheses are considered in this example:
\begin{equation}
  \mathcal{H}_1 \colon \bm{\beta} = \bm{\beta}_1 = (0.75, 0.25), \qquad \mathcal{H}_2 \colon \bm{\beta} = \bm{\beta}_2 = (0.5, 0.5).
  \label{eq:hyp_example2}
\end{equation}
Note that for $\bm{\beta}_2 = (0.5, 0.5)$ the binomial AR(1) process reduces to a process of independent binomial random variables with distribution $B(M,0.5)$. Hence, the two hypotheses in \eqref{eq:hyp_example2} correspond to a test for dependencies in the observed data. The aim is to solve the Kiefer--Weiss problem for the hypotheses in \eqref{eq:hyp_example2}, that is, to design a sequential test whose worst-case expected run length over all possible random processes is minimal. Consequently, the uncertainty sets for the conditional distributions $P_{\theta}$ are chosen as $\mathcal{P}_{\theta} = \mathcal{M}_{\mu}$ for all $\theta \in \Omega_\theta$. Note that this type of uncertainty can be interpreted as a special case of the density band model in \eqref{eq:band_model}, with $p' = 0$ and $p'' = 1$, or as an outlier model with contamination ratio $\varepsilon = 1$. In analogy to the previous example, the reference measure is set to $\mu = P_{\theta}^{(2)} = B(7, 0.5)$, so that $z_2 = 1$ and $\rhofunc$ becomes a function of $(z_0, z_1)$ only. The initial state of the sufficient statistic is set to $\theta_0 = 3$.

Since the hypotheses in \eqref{eq:hyp_example2} are simple, a regular sequential probability ratio test with log-likelihood ratio thresholds $B < 0 < A$ can be applied as well. However, under the above uncertainty model, its worst-case expected run length is infinite. In order to see this, consider a deterministic process that alternates between two observations, $x_{(1)}$ and $x_{(2)}$, which are chosen such that their log-likelihood ratios satisfy
\begin{equation}
  0 < \log \frac{P_{x_{(1)}, \bm{\beta}_1}(x_{(2)})}{P_{x_{(1)}, \bm{\beta}_0}(x_{(2)})} < A \quad \text{and} \quad B < \log \frac{P_{x_{(2)}, \bm{\beta}_1}(x_{(1)})}{P_{x_{(2)}, \bm{\beta}_0}(x_{(1)})} < 0.
\end{equation}
For this process the log-likelihood ratio increments keep canceling each other out so that neither of the thresholds is ever crossed. A minimax robust test, by contrast, makes it possible to leverage the increased efficiency of sequential tests while at the same time having a bounded worst-case run length.

\begin{figure}[tb]
  \centering
  \includegraphics[width=0.5\textwidth]{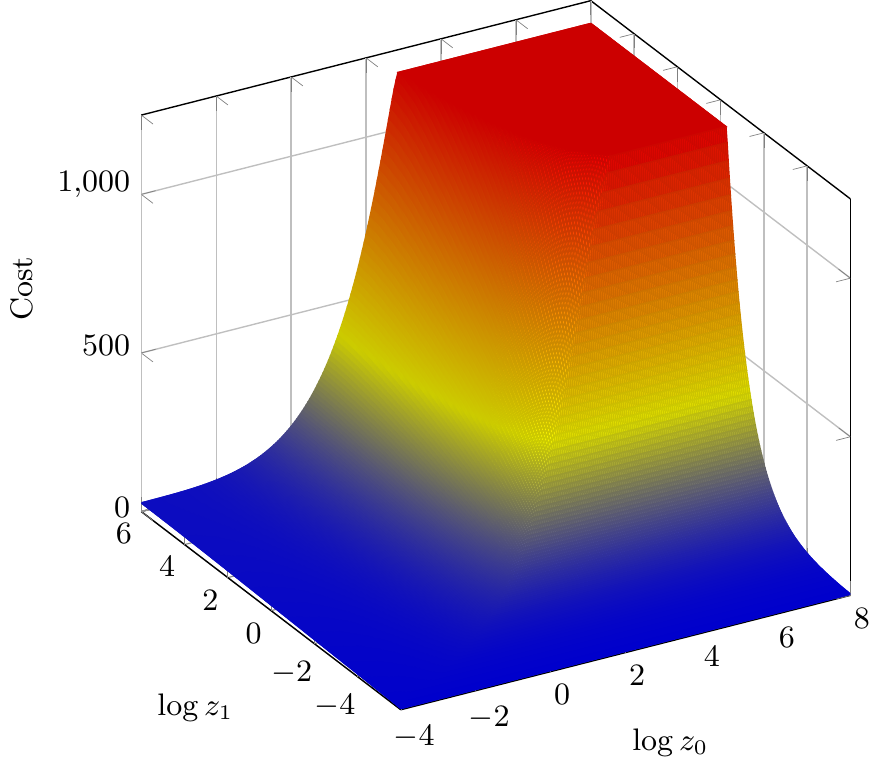}%
  \includegraphics[width=0.5\textwidth]{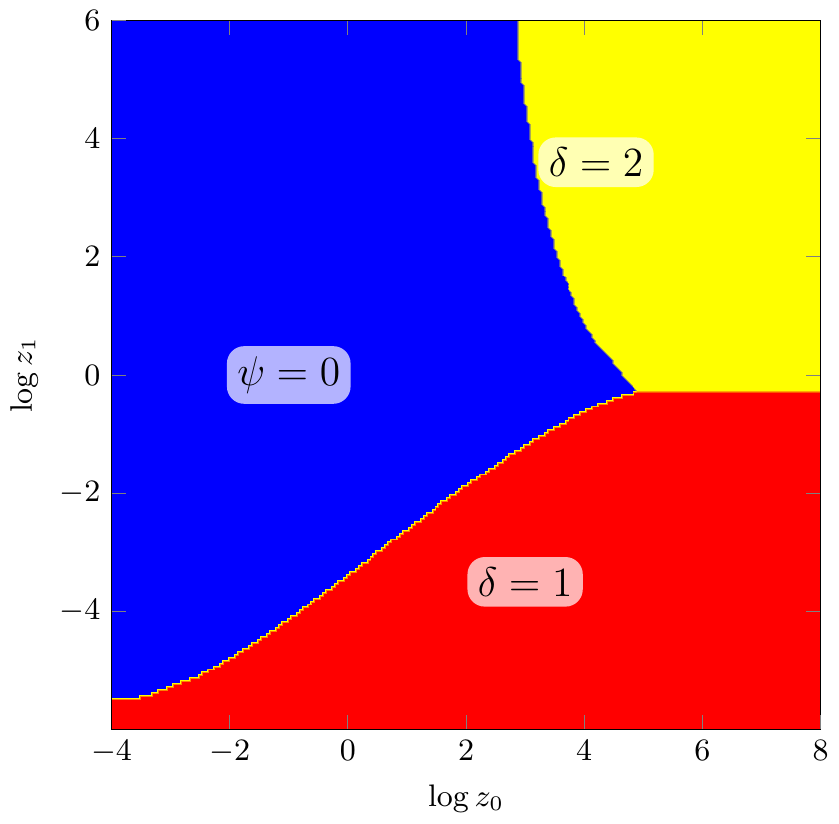}
  \caption{Segment of the optimal cost function $\rhofunc(\z, 3)$ in logarithmic scale and the corresponding testing policy as a function of the log-likelihood ratios.}
  \label{fig:opt_policy_2}
\end{figure}

In order to obtain a numerical solution, the likelihood ratio plane $(z_0, z_1)$ was discretized on $[-4, 8] \times [-6, 6]$ using $241 \times 241$ uniformly spaced grid points. The iterative design procedure converged after four iterations. The optimal weights were found to be $\bm{\lambda}^* \approx (1526.38, 1178.24)$. The resulting cost function $\rho_{\bm{\lambda}^*}$, as well as the corresponding testing policy, are depicted in Figure~\ref{fig:opt_policy_2}. Both are distinctly different from their counterparts in the first example. First, the stopping region is no longer a corridor but is of a conic shape with the thresholds tightening as $\log z_0$ increases. Second, it is noteworthy that for $\log z_0$ greater than approximately five, the sequential test reduces to a single threshold test, which is an indicator for the test being truncated under certain conditions and is in line with the goal of minimizing the worst-case expected run length. For different values of $\theta$ slight changes in the location of the decision regions can be observed but the overall shape remains the same.

\begin{figure}[tb]
  \centering
  \includegraphics[width=0.5\textwidth]{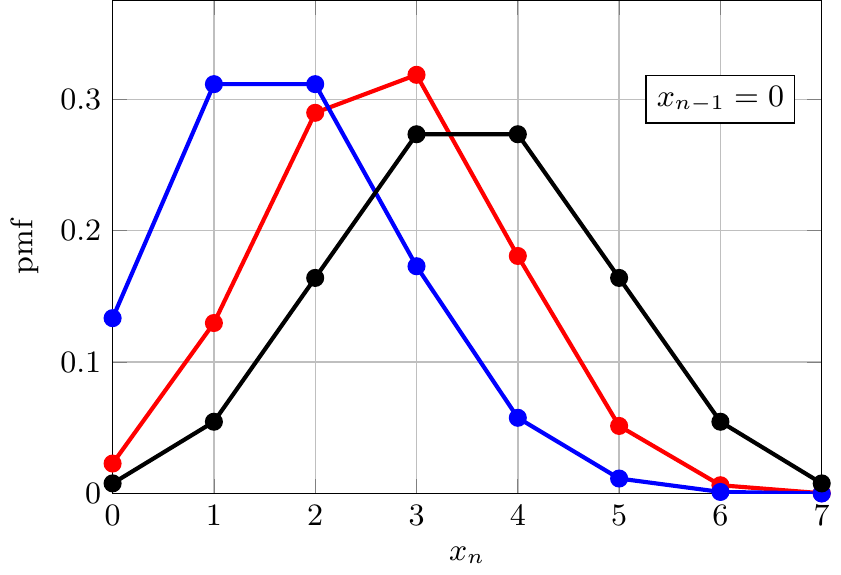}%
  \includegraphics[width=0.5\textwidth]{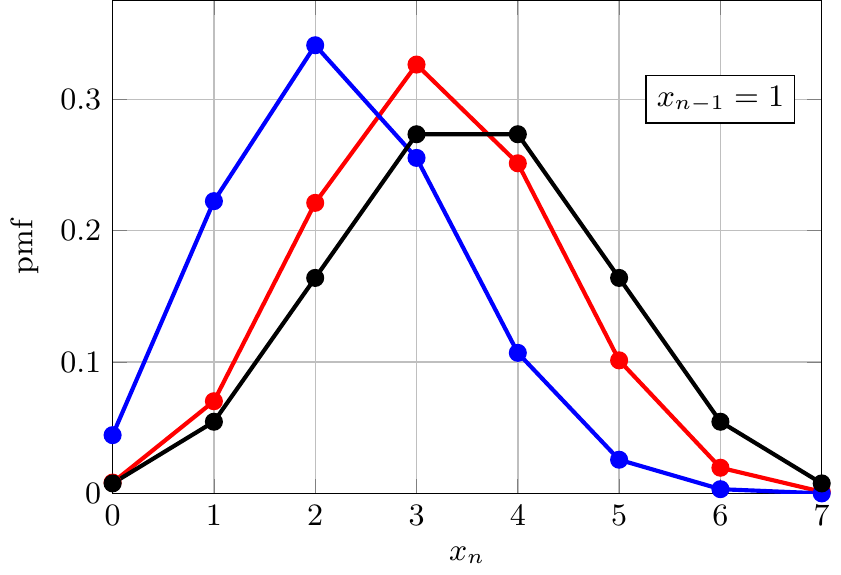}\\
  \includegraphics[width=0.5\textwidth]{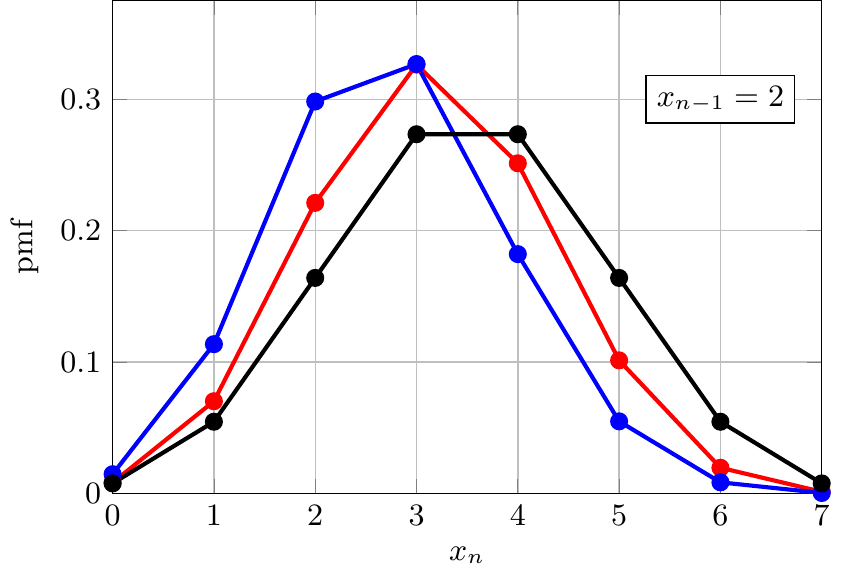}%
  \includegraphics[width=0.5\textwidth]{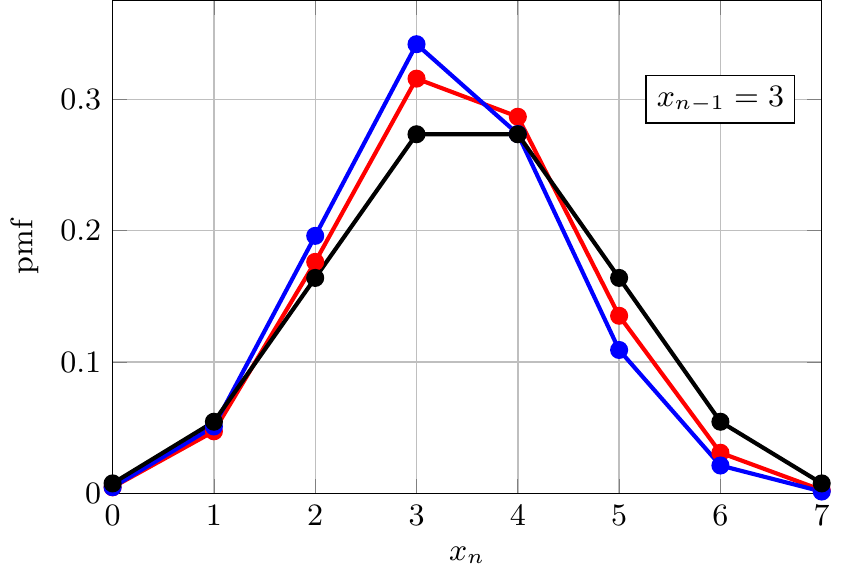}\\
  \includegraphics[width=0.5\textwidth]{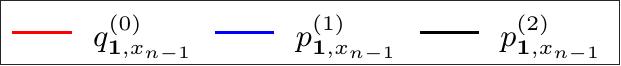}
  \caption{Examples of least favorable distributions in the state $(\bm{1}, \theta)$ for different values of the sufficient statistic $\theta$.}
  \label{fig:lfds_2}
\end{figure}

Four examples of conditional distributions that are least favorable with respect to the expected run length are depicted in Figure~\ref{fig:lfds_2}. Note that a linear interpolation is used to connect the point masses; although this a potentially misleading representation, it helps to make the differences in shape more recognizable. In the depicted examples the least favorable distributions are conditioned on $\z = \bm{1}$ and on different values for $\theta$, which corresponds to the previous observation $x_{n-1}$. It can be observed how the least favorable distribution, depicted in red, adapts to the state of the test statistic in such a way that it is equally similar to both hypotheses. Moreover, it is noteworthy that, even without bounds on the densities, the least favorable distributions do not reduce to a single point mass, meaning that for any given state of the test statistic there is no single least favorable observation. This is a consequence of the fact that the least favorable distributions under density band uncertainty turn out to be equalizers with respect to the respective performance measure, so that in this case all observations lead to the same expected run length---see \cite{Fauss2017_isi} for a more detailed discussion.

In analogy to the first example, $\num{10000}$ Monte Carlo simulations were performed using the testing policy depicted in Figure~\ref{fig:opt_policy_2} in order to verify the numerical results. The confusion matrix as well the average run length of the tests are shown in Table~\ref{table:mc_results_2}. The average run length under the least favorable distribution $\mathbb{Q}_0$ was obtained to be approximately $74.66$ samples. Interestingly, a fixed sample size test between $\mathcal{H}_1$ and $\mathcal{H}_2$ using $75$ samples results in error probabilities close to $1\%$ as well. This raises the question whether there is a relation between the worst-case expected run length of a minimax sequential test and the number of samples required by the equivalent fixed-sample size test. More precisely, since the latter is an upper bound on the former, the question is whether the least favorable distributions attain this bound; this conjecture will be investigated in future work. Regardless of whether or not the conjecture holds, the minimax sequential test still achieves significant reductions in the average sample size under both hypotheses, in particular $\mathcal{H}_2$, as is clear by inspection of Table~\ref{table:mc_results_2}.

\begin{table}[tb]
  \centering
  \includegraphics{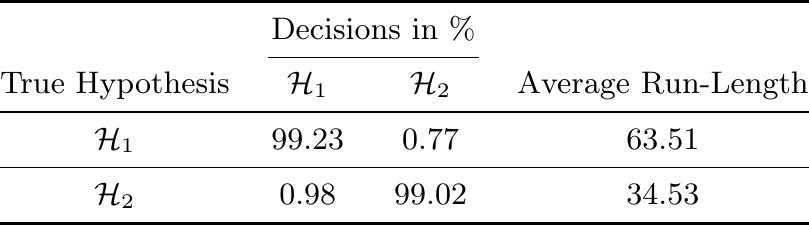}%
  \caption{Results of $10\,000$ Monte Carlo runs of the minimax optimal test using the policy depicted in Figure\ref{fig:opt_policy_2}. The target detection probability is $99\%$.}
  \label{table:mc_results_2}
\end{table}

\appendix

\section{Proof of Theorem~\titleref{th:rho_sum}}
\label{apd:proof_rho_sum}

In order to prove Theorem~\ref{th:rho_sum}, it suffices to show that the right hand side of \eqref{eq:rho_sum} solves the integral equation \eqref{eq:optimal_rho}. Since $\rhofunc$ is unique, this implies that both functions are identical. 

Two auxiliary results are used in the proof:
\begin{enumerate}
  \item[(I)] For every function of the form $z_k f(\z,\theta)$, with $f \colon \rhodomain \to \Rnng$, it holds that
             \begin{align*}
                \int z_k f(\z,\theta) \, \dint \mu_{\z,\theta} &= \int z_k \, p_{\theta}^{(k)}(x) \, f\left(\z \p_{\theta}(x), \xi_{\theta}(x) \right) \, \mu(\dint x) \\
                &= z_k \int f\left(\z \p_{\theta}(x), \xi_{\theta}(x) \right) \, P_\theta^{(k)}(\dint x) \\
                &= z_k \int f \, \dint P_{\z,\theta}^{(k)}.
             \end{align*}
        
  \item[(II)] For $\pi = (\psi, \bm{\delta}) \in \PiLambda^*(\pmb{\Pbb})$, it follows from Corollary~\ref{col:optimal_policies} that
               \begin{align*}
                 \sum_{k=1}^K \lambda_k z_k (1-\delta_k) = \sum_{k=1}^K \lambda_k z_k - \max_{k=1,\ldots,K}\{ \lambda_k z_k \} = \gfunc.
               \end{align*}
\end{enumerate}
Using these properties and the Chapman--Kolmogorov equations \eqref{eq:ck_runlength} and \eqref{eq:ck_errors}, it follows that for all $(\psi, \bm{\delta}) \in \PiLambda^*(\pmb{\Pbb})$
\begin{align*}
  \rhofunc &= z_0 \gamma_{\pi,\Pbb_0} + \sum_{k=1}^K \lambda_k z_k \alpha_{\pi,\Pbb_k}^{(k)} \\
  &= z_0 (1-\psi)  \left( 1 + \int \gamma_{\pi,\Pbb_0} \,\dint P_{\z,\theta}^{(0)} \right) 
    + \sum_{k=1}^K \lambda_k z_k \left( \psi (1-\delta_k) + (1-\psi) \int \alpha_{\pi,\Pbb_k}^{(k)} \,\dint P_{\z,\theta}^{(k)} \right) \\
  &\refeqq{(I)} (1-\psi) \left( z_0 + \int z_0 \gamma_{\pi,\Pbb_0} \,\dint \mu_{\z,\theta} \right) 
    + \sum_{k=1}^K \lambda_k z_k \psi (1-\delta_k) + \sum_{k=1}^K (1-\psi) \int \lambda_k z_k \alpha_{\pi,\Pbb_k}^{(k)} \,\dint \mu_{\z,\theta} \\
  &\refeqq{(II)} \psi \gfunc + (1-\psi) \left( z_0 + \int \left( z_0 \gamma_{\pi,\Pbb_0} + \sum_{k=1}^K \lambda_k z_k \alpha_{\pi,\Pbb_k}^{(k)} \right) \,\dint \mu_{\z,\theta} \right) \\
  &= \psi \gfunc + (1-\psi) \left( z_0 + \int \rhofunc \,\dint \mu_{\z,\theta} \right) \\
  &= \min \left\{ \gfunc \,,\, z_0 + \int \rhofunc \,\dint \mu_{\z, \theta} \right\},
\end{align*}
where the last equality follows again from Corollary~\ref{col:optimal_policies}. This concludes the proof.

\section{Proof of Lemma~\titleref{lm:rho_properties}}
\label{apd:proof_rho_properties}

Consider the sequence of functions $\bigl( \rhofunc^n \bigr)_{n \geq 0}$ with $\rhofunc^0 = \gfunc$ and
\begin{equation}
  \rhofunc^n = \min \; \left\{ \gfunc(\z) \; , \; z_0 + \int \rhofunc^{n-1} \, \dint \mu_{\z,\theta} \right\}.
  \label{eq:rho_sequence}
\end{equation}
Assume that $\rhofunc^n \leq \rhofunc^{n-1}$ for some $n \geq 0$. It then holds that
\begin{align}
  \rhofunc^{n+1} &= \min \; \left\{ \gfunc(\z) \; , \; z_0 + \int \rhofunc^{n} \, \dint \mu_{\z,\theta} \right\} \\
  &\leq \min \; \left\{ \gfunc(\z) \; , \; z_0 + \int \rhofunc^{n-1} \, \dint \mu_{\z,\theta} \right\} = \rhofunc^{n}
\end{align}
so that $\bigl( \rhofunc^n \bigr)_{n \geq 0}$ is non-increasing. Moreover, the sequence is non-negative and hence bounded from below. This implies that it converges pointwise to a unique limit. Since this result is used repeatedly in the paper, it is fixed in the next lemma.
\begin{lemma}
  Let $(\Omega, \Filter)$ be a measurable space and let $\bigl( f_n \bigr)_{n \geq 0}$, with $f_n \colon \Omega \to \mathbb{R}$, be a sequence of functions. If this sequence is non-increasing, that is, $f_n \leq f_{n-1}$ for all $n \geq 0$, and a function $g > -\infty$ exists such that $f_n \geq g$ for all $n \geq 0$, then the pointwise limit
  \begin{equation}
    \lim_{n \to \infty} \; f_n = f
  \end{equation}
  exists and is unique. The same result holds if the sequence is non-decreasing, that is, $f_n \geq f_{n-1}$ for all $n \geq 0$, and a function $g < \infty$ exists such that $f_n \leq g$.
  \label{lm:monotonic_convergence}
\end{lemma}
Lemma~\ref{lm:monotonic_convergence} is a well-known result in Real Analysis. It follows from the fact that for every $\omega \in \Omega$, it holds that $\bigl( f_n(\omega) \bigr)_{n \geq 0}$ is a monotonic and bounded sequence of real numbers. Sequences of this type are guaranteed to converge to a unique limit; see, for example, \cite[Theorem 3.14]{Rudin1976} and \cite[Theorem 3.1.4]{Mukherjee2005}. This immediately implies the statement in Lemma~\ref{lm:monotonic_convergence}.

From Lemma~\ref{lm:monotonic_convergence} it follows that the sequence $\bigl( \rhofunc^n \bigr)_{n \geq 0}$ converges pointwise to a unique limit $\rhofunc$ for $n \to \infty$; also compare \cite[Lemma 4, Lemma 5]{Novikov2009_multiple_hypotheses} and \cite[Appendix A]{Fauss2015}. By definition \eqref{eq:g_function}, $\gfunc$ is non-decreasing, concave, and homogeneous of degree one in $\z$. Using this as a basis, it can be shown via induction that these properties carry over to $\rhofunc$. Here, only the concavity of $\rhofunc$ is proven in detail; the other properties can be shown analogously. 

Assume that $\rhofunc^n$ is concave for some $n \geq 0$, that is,
\begin{equation}
  \rhofunc^n\bigl(\kappa \z' + (1-\kappa) \z,\theta \bigr) \geq \kappa \rhofunc^n(\z',\theta) + (1-\kappa)\rhofunc^n(\z, \theta)
  \label{eq:concavity_rho_n}
\end{equation}
for all $\z, \z' \in \rhodomain$ and all $\kappa \in [0,1]$. From the fact that $\gfunc$ and the minimum are concave functions it follows that
\begin{align*}
  \rhofunc^{n+1}(\kappa \z' + (1-\kappa) \z,\theta) &= \min \left\{ \gfunc(\kappa \z' + (1-\kappa) \z) \; , \; z_0 + \int \rhofunc^n \dint \mu_{\kappa \z' + (1-\kappa) \z, \theta} \right\} \\
    &\geq \min \biggl\{ \kappa \gfunc(\z') + (1-\kappa) \gfunc(\z) \,,\, \\
    &\qquad z_0 + \kappa \! \int \rhofunc^n \dint \mu_{\z', \theta} + (1-\kappa) \! \int \rhofunc^n \dint \mu_{\z, \theta} \biggr\} \\
    &\geq \kappa \min \left\{ \gfunc(\z') \,,\, z_0 + \int \rhofunc^n \dint \mu_{\z',\theta} \right\} \\ 
    &\qquad + (1-\kappa) \min \left\{ \gfunc(\z) \,,\, z_0 +  \int \rhofunc^n \dint \mu_{\z, \theta} \right\} \\
    &= \kappa \rhofunc^{n+1}(\z',\theta) + (1-\kappa) \rhofunc^{n+1}(\z, \theta)
\end{align*}
so that $\rhofunc^{n+1}$ is concave as well. Taking the limit on both sides of \eqref{eq:concavity_rho_n} yields
\begin{equation}
  \rhofunc\bigl(\kappa \z' + (1-\kappa) \z,\theta \bigr) \geq \kappa \rhofunc(\z',\theta) + (1-\kappa)\rhofunc(\z, \theta).
\end{equation}
This concludes the proof.

\section{Proof of Theorem~\titleref{th:rho_derivatives}}
\label{apd:proof_rho_derivatives}

Since $\rhofunc$ is concave, its generalized partial derivatives in the sense of \eqref{eq:partial_diff} exist. Moreover, for $\pi \in \PiLambda^*(\pmb{\Pbb})$, $\rhofunc$ can be written as (compare Appendix~\ref{apd:proof_rho_sum})
\begin{equation}
  \rhofunc = \min \{ \gfunc \,,\, \dfunc \,\} = \psi \left( \sum_{k+1}^K \lambda_k z_k (1-\delta_k) \right) + (1-\psi) \dfunc
  \label{eq:min_g_d}
\end{equation}
where $\dfunc \colon \rhodomain \to \Rnng$ is defined as
\begin{equation}
  \dfunc(\z,\theta) \coloneqq z_0 + \int \rhofunc \,\dint \mu_{\z,\theta}.
\end{equation}
Exploiting the coupling between $\rhofunc$ and $\dfunc$ is key to the proof of Theorem~\ref{th:rho_derivatives}. The argument used here is based on a generalized version of Leibniz's integral rule, which is given in the next lemma.
\begin{lemma}[Generalized Leibniz integral rule]
  \label{lm:leibniz}
  Let $(\Omega,\mathcal{F})$ be a measurable space and let $f\colon \mathbb{R}^K \times \Omega \to \mathbb{R}$ be a convex or concave function. If $f(\bm{y},\omega)$ is $\mu$-integrable for all $\bm{y} \in \mathbb{R}^K$, it holds that
  \begin{equation*}
    \partial_{y_k} \left( \int_{\Omega} f(\bm{y}, \omega) \, \mu(\dint \omega) \right) = \int_{\Omega} \partial_{y_k} f(\bm{y}, \omega) \, \mu(\mathrm{d} \omega),
  \end{equation*}
  where the integral on the right hand side is a short-hand notation for the set of integrals over all feasible partial derivatives of $f$, i,e.,
  \begin{equation*}
    \int_{\Omega} \partial_{y_k} f(\bm{y}, \omega) \, \mu(\dint \omega) \coloneqq \left\{ c \in \mathbb{R} : \exists f_{y_k} \in \partial_{y_k} f : c = \int_{\Omega} f_{y_k}(\bm{y}, \omega) \, \mu(\mathrm{d} \omega)  \right\}.
  \end{equation*}
\end{lemma}
The generalized Leibniz integral rule is proven, for example, in \cite[Theorem 23]{Rockafellar1974}. Extensions and variations are given in
\cite{Papageorgiou1997} and \cite{Chieu2009}.

Since
\begin{equation}
  \int \rhofunc \, \dint \mu_{\z,\theta} \leq \int \gfunc \, \dint \mu_{\z,\theta} \leq \sum_{k=1}^K \lambda_k z_k < \infty,
\end{equation}
$\rhofunc$ is $\mu_{\z,\theta}$-integrable for all $(\z,\theta) \in \rhodomain$. Hence, Leibniz's integral rule applies to $\dfunc$ so that
\begin{align*}
  \partial_{z_k} \dfunc(\z,\theta) &= \partial_{z_k} \left( z_0 + \int \rhofunc \dint \mu_{\z, \theta} \right) \\
  &= \partial_{z_k} z_0 + \int \partial_{z_k} \Bigl( \rhofunc \left( \z \p_{\theta}(x), \xi_{\theta}(x) \right) \Bigr) \, \mu(\dint x) \\
  &= \partial_{z_k} z_0 + \int \partial_{z_k} \rhofunc \left(\z \p_{\theta}(x), \xi_{\theta}(x) \right) p_{\theta}^{(k)}(x) \, \mu(\dint x) \\
  &= \partial_{z_k} z_0 + \int \partial_{z_k} \rhofunc \left(\z \p_{\theta}(x), \xi_{\theta}(x) \right) \, P_{\theta}^{(k)}(\dint x) \\
  &= \partial_{z_k} z_0 + \int \partial_{z_k} \rhofunc \, \dint P_{\z, \theta}^{(k)}
\end{align*}
By expressing $\rhofunc$ as in \eqref{eq:min_g_d} and taking the partial generalized derivatives, the following set-valued integral equations \cite{Berenguer2015} are obtained:
\begin{equation}
  \partial_{z_0} \rhofunc(\z,\theta) = (1-\psi) \left( 1 + \int \partial_{z_0} \rhofunc \, \dint P_{\z, \theta}^{(0)} \right)
  \label{eq:set_integral_diff_0}
\end{equation}
for $k = 0$ and
\begin{equation}
  \partial_{z_k} \rhofunc(\z,\theta) = \psi \, \lambda_k (1-\delta_k) + (1-\psi) \int \partial_{z_k} \rhofunc \, \dint P_{\z, \theta}^{(k)}
  \label{eq:set_integral_diff_k}
\end{equation}
for $k = 1,\ldots, K$. Read from ``right to left'', \eqref{eq:set_integral_diff_0} and \eqref{eq:set_integral_diff_k} state that inserting any function $r'_k \in \partial_{z_k} \rhofunc$  into the integral on the right hand side yields another function $r_k \in \partial_{z_k} \rhofunc$ on the left hand side. Read from ``left to right'', \eqref{eq:set_integral_diff_0} and \eqref{eq:set_integral_diff_k} state that given any function $r_k \in \partial_{z_k} \rhofunc$ on the left hand side, a function $r'_k \in \partial_{z_k} \rhofunc$ exists such that the right hand side evaluates to $r_k$.

The above characterization of the generalized differentials, which follows solely from the concavity and integrability of $\rhofunc$, already implies both statements in Theorem~\ref{th:rho_derivatives}. By inspection of \eqref{eq:ck_runlength_compact} and \eqref{eq:ck_errors_compact}, it can be seen that $\gamma_{\pi,\Pbb_0}$ and $\lambda_k \alpha_{\pi,\Pbb_k}^{(k)}$ are solutions of \eqref{eq:set_integral_diff_0} and \eqref{eq:set_integral_diff_k}, respectively, for all $\pi \in \PiLambda^*(\pmb{\Pbb})$. This yields the first statement in Theorem~\ref{th:rho_derivatives}.

The second part of Theorem~\ref{th:rho_derivatives} is proven by showing that the sets in the statement are subsets of each other. The details are only given for $k \geq 1$ since the proof for $k = 0$ follows analogously. From the above results it is clear that $\lambda_k \alpha_{\pi,\Pbb_k}^{(k)} \in \partial_{z_k} \rhofunc$ for all $\pi \in \PiLambda^*(\pmb{\Pbb})$. By definition, this implies that $\lambda_k \alpha_{\pi,\Pbb_k}^{(k)}(\z,\theta) \in \partial_{z_k} \rhofunc(\z,\theta)$ for all $(\z,\theta) \in \rhodomain$ and all $\pi \in \PiLambda^*(\pmb{\Pbb})$, that is,
\begin{equation}
  \Bigl\{\, \lambda_k \alpha_{\pi,\Pbb_k}^{(k)}(\z,\theta) : \pi \in \PiLambda^*(\pmb{\Pbb}) \,\Bigr\} \subset \partial_{z_k} \rhofunc(\z,\theta).
\end{equation}
In order to show the converse, the following lemma is useful.
\begin{lemma}
  Given two policies $\pi, \pi' \in \PiLambda^*(\pmb{\Pbb})$ that satisfy
  \begin{equation}
    a = \alpha_{\pi,\Pbb_k}^{(k)}(\z,\theta) \leq \alpha_{\pi', \Pbb_k}^{(k)}(\z,\theta) = a'
  \end{equation}
  for some $(\z,\theta) \in \rhodomain$, it holds that for every $\tilde{a} \in [a, a']$ there exits a policy $\tilde{\pi} \in \PiLambda^*(\pmb{\Pbb})$ such that
  \begin{equation}
    \alpha_{\tilde{\pi}, \Pbb_k}^{(k)}(\z,\theta) = \tilde{a}.
  \end{equation}
  \label{lm:mixed_policies}
\end{lemma}
The lemma can be shown by considering a randomized policy $\tilde{\pi}$, which, at each time instant, is chosen to be $\pi$ with probability $\kappa$ and $\pi'$ with probability $(1-\kappa)$, where $\kappa \in [0,1]$. The conditional probability of erroneously deciding against $\mathcal{H}_k$ when using this mixed policy is given by the integral equation
\begin{equation}
  \alpha_{\tilde{\pi}, \Pbb_k}^{(k)} = \kappa \psi (1-\delta_k) + (1-\kappa)\psi' (1-\delta_k') + \int \alpha_{\tilde{\pi}, \Pbb_k}^{(k)} \dint P_{\z, \theta}^{(k)},
\end{equation}
which has the unique solution $\alpha_{\tilde{\pi}, \Pbb_k}^{(k)} = \kappa \alpha_{\pi,\Pbb_k}^{(k)} + (1-\kappa) \alpha_{\pi', \Pbb_k}^{(k)}$ so that 
\begin{equation}
  \alpha_{\tilde{\pi}, \Pbb_k}^{(k)}(\z,\theta) = \kappa \alpha_{\pi,\Pbb_k}^{(k)}(\z,\theta) + (1-\kappa) \alpha_{\pi', \Pbb_k}^{(k)}(\z,\theta) = \kappa a + (1-\kappa) a'.
\end{equation}
The Lemma follows.

Using Lemma~\ref{lm:mixed_policies}, the desired result can be shown by contradiction. Assume that $r_k \in \partial_{z_z} \rhofunc$ exists such that for some $(\z, \theta)$ it holds that
\begin{equation}
  r_k(\z, \theta) \notin \Bigl\{\, \lambda_k \alpha_{\pi,\Pbb_k}^{(k)}(\z,\theta) : \pi \in \PiLambda^*(\pmb{\Pbb}) \,\Bigr\}.
  \label{eq:r_assumption}
\end{equation}
By \eqref{eq:set_integral_diff_k}, a policy $\pi' \in \PiLambda^*(\pmb{\Pbb})$ and a function $r'_k \in \partial_{z_k} \rhofunc$ exists such that
\begin{equation}
  r_k = \psi' \, \lambda_k (1-\delta'_k) + (1-\psi') \int r'_k \, \dint P_{\z,\theta}^{(k)}.
\end{equation}
Now consider the sequence of functions defined by
\begin{equation}
  \overline{r}_k^{n} = \sup_{\pi \in \PiLambda^*(\pmb{\Pbb})} \; \left\{ \psi \, \lambda_k (1-\delta_k) + (1-\psi) \int \overline{r}_k^{n-1} \, \dint P_{\z,\theta}^{(k)} \right\}
  \label{eq:r_sequence}
\end{equation}
with $\overline{r}_k^{0} = r_k'$. Assume that $\overline{r}_k^{n} \geq \overline{r}_k^{n-1}$ for some $n \geq 0$. Via induction, it follows that
\begin{align}
  \overline{r}_k^{n+1} &= \sup_{\pi \in \PiLambda^*(\pmb{\Pbb})} \; \left\{ \psi \, \lambda_k (1-\delta_k) + (1-\psi) \int \overline{r}_k^{n} \, \dint P_{\z,\theta}^{(k)} \right\} \\
  \geq &= \sup_{\pi \in \PiLambda^*(\pmb{\Pbb})} \; \left\{ \psi \, \lambda_k (1-\delta_k) + (1-\psi) \int \overline{r}_k^{n-1} \, \dint P_{\z,\theta}^{(k)} \right\} = \overline{r}_k^n.
\end{align}
Since the induction basis is satisfied by construction of $\overline{r}_k^{0}$, it follows that $\bigl( \overline{r}_k^{n} \bigr)_{n \geq 0}$ is a non-decreasing sequence. Moreover, since $\rhofunc$ is concave in $\z$, it holds that
\begin{align}
  \sup \, \partial_{z_k} \rhofunc(\z, \theta) \leq \sup \, \partial_{z_k} \rhofunc(\bm{0}, \theta) = \partial_{z_k} \gfunc(\z) = \lambda_k
\end{align}
so that $\overline{r}_k^{0} = r_k'$ is bounded from above by $\lambda_k$. Again, using this as an induction basis, it follows that
\begin{align}
  \overline{r}_k^{n} &= \sup_{\pi \in \PiLambda^*(\pmb{\Pbb})} \; \left\{ \psi \, \lambda_k (1-\delta_k) + (1-\psi) \int \overline{r}_k^{n-1} \, \dint P_{\z,\theta}^{(k)} \right\} \\
  &\leq \sup_{\pi \in \PiLambda^*(\pmb{\Pbb})} \; \left\{ \psi \, \lambda_k (1-\delta_k) + (1-\psi) \lambda_k \right\} \leq \lambda_k.
\end{align}
so that $\bigl( \overline{r}_k^{n} \bigr)_{n \geq 0}$ is bounded by $\lambda_k$. Hence, by Lemma~\ref{lm:monotonic_convergence}, the limit $\lim_{n \to \infty} \overline{r}_k^{n} = \overline{r}_k$ exist and it solves the integral equation
\begin{equation}
  \overline{r}_k = \sup_{\pi \in \PiLambda^*(\pmb{\Pbb})} \; \left\{ \psi \, \lambda_k (1-\delta_k) + (1-\psi) \int \overline{r}_k \, \dint P_{\z,\theta}^{(k)} \right\}.
  \label{eq:integral_limit_r}
\end{equation}
By inspection of $\PiLambda^*(\pmb{\Pbb})$ in Corollary~\ref{col:optimal_policies}, the supremum is achieved by the decision and stopping rules defined by the decision functions 
\begin{equation}
    \delta_k(\z) = \mathcal{I}\left(\left\{\Gfunc(\z) = \gfunc(\z)\right\}\right)
\end{equation}
and
\begin{equation}
  \overline{\psi}(\z, \theta) = 
    \begin{dcases}
      \mathcal{I}\left(\{\gfunc(\z) < \rhofunc(\z, \theta)\}\right), & \lambda_k (1-\overline{\delta}_k(\z)) < \int \overline{r}_k \, \dint P_{\z,\theta}^{(k)} \\
      \mathcal{I}\left(\{\gfunc(\z) \leq \rhofunc(\z, \theta)\}\right), & \lambda_k (1-\overline{\delta}_k(\z)) \geq \int \overline{r}_k \, \dint P_{\z,\theta}^{(k)}
    \end{dcases}
\end{equation}
so that \eqref{eq:integral_limit_r} can equivalently be written as
\begin{equation}
  \overline{r}_k = \overline{\psi} \, \lambda_k (1-\overline{\delta}_k) + (1-\overline{\psi}) \int \overline{r}_k \, \dint P_{\z,\theta}^{(k)}.
\end{equation}
Consequently, it holds that $\overline{r}_k = \lambda_k \alpha_{\overline{\pi},\mathbb{P}_k}^{(k)} \geq r_k$. Analogously, by replacing the supremum in \eqref{eq:r_sequence} with the infimum, a sequence $\bigl( \underline{r}_k^{n} \bigr)_{n \geq 0}$ can be constructed that is non-increasing, bounded from below by zero and, therefore, converges to a unique limit $\underline{r}_k = \lambda_k \alpha_{\underline{\pi},\mathbb{P}_k}^{(k)} \leq r_k$. This implies that
\begin{equation}
  \lambda_k \alpha_{\underline{\pi},\mathbb{P}_k}^{(k)}(\z, \theta) \leq r_k(\z, \theta) \leq \lambda_k \alpha_{\overline{\pi},\mathbb{P}_k}^{(k)}(\z, \theta).
\end{equation}
From Lemma~\ref{lm:mixed_policies} it now follows that a policy $\tilde{\pi} \in \PiLambda^*(\pmb{\Pbb})$ exists such that 
\begin{equation}
  \lambda_k \alpha_{\tilde{\pi},\mathbb{P}_k}^{(k)} = r_k(\z, \theta),
\end{equation}
which contradicts the assumption \eqref{eq:r_assumption}. Hence, it holds that
\begin{equation}
  \left\{\, \alpha_{\pi,\Pbb_k}^{(k)}(\z,\theta) : \pi \in \PiLambda^*(\pmb{\Pbb}) \,\right\} \supset \partial_{z_k} \rhofunc(\z,\theta),
\end{equation}
for all $(\z,\theta) \in \rhodomain$. This concludes the proof of the second statement in Theorem~\ref{th:rho_derivatives}.

\section{Proof of Theorem~\titleref{th:optimal_constrained_test}}
\label{apd:proof_optimal_constrained_test}

The outline of the proof is as follows: First, it is shown that if $\bm{\lambda}^*$ satisfies \eqref{eq:lambda_opt}, it holds that $\gamma(\pi^*,\pmb{\Pbb}) \in \partial_{z_0} \rho_{\bm{\lambda}^*}(\bm{1}, \theta_0)$ and that $\overline{\alpha}_k \in \partial_{z_k} \rho_{\bm{\lambda}^*}(\bm{1}, \theta_0)$. Theorem~\ref{th:optimal_constrained_test} then follows as a consequence of Theorem~\ref{th:rho_derivatives}.

In order to show the first part, the following lemma is useful.
\begin{lemma}
  For $\rhofunc$ defined in Theorem~\ref{th:optimal_cost_function} it holds that for all $(\z, \theta) \in  \rhodomain$ and all $k = 0, \ldots, K$
  \begin{equation}
    \lambda_k \partial_{\lambda_k} \rhofunc(\z, \theta) = z_k \partial_{z_k} \rhofunc(\z, \theta)
  \end{equation}
  where $\lambda_0 = 1$ is defined to unify notation.
  \label{lm:diff_lambda_z}
\end{lemma}
In order to show the lemma, assume that a function $\tilde{\rho} \colon \rhodomain \to \Rnng$ exists such that for all $(\z,\theta) \in \rhodomain$ it holds that
\begin{equation}
  \rhofunc(\z,\theta) = \tilde{\rho}(\bm{\lambda} \z, \theta), 
  \label{def:rho_tilde}
\end{equation}
Given that \eqref{def:rho_tilde} holds, it immediately follows that
\begin{align*}
  z_k \, \partial_{z_k} \rhofunc(z,\theta) &= z_k \, \partial_{z_k} \tilde{\rho}(\bm{\lambda} \z, \theta) \\
  &= \lambda_k z_k \, \partial_{\lambda_k z_k} \tilde{\rho}(\bm{\lambda} \z, \theta) \\
  &= \lambda_k \, \partial_{\lambda_k} \tilde{\rho}(\bm{\lambda} \z, \theta) = \lambda_k \partial_{\lambda_k} \rhofunc(\z,\theta),
\end{align*}
where $\partial_{\lambda_k z_k} \tilde{\rho}$ denotes the generalized differential of $\tilde \rho$ with respect to ${\lambda_k z_k}$. The existence of $\tilde{\rho}$ can be shown via induction. Let the sequence $\rhofunc^n$ be as defined in \eqref{eq:rho_sequence} and assume that \eqref{def:rho_tilde} holds for some $n \geq 0$, that is, a function $\tilde{\rho}^n$ exists such that $\rhofunc^n(z,\theta) = \tilde{\rho}^n(\bm{\lambda} \z,\theta)$. It then follows that
\begin{align*}
  \rhofunc^{n+1}(\z,\theta) &= \min \left\{ \gfunc(\z) \; , \; 1 + \int \rhofunc^n \left( \z \p_\theta(x), \xi_{\theta_n}(x) \right) \dint \mu(x) \right\} \\
  &= \min \left\{ \tilde{g}(\bm{\lambda} \z) \; , \; 1 + \int \tilde{\rho}^n \left( \bm{\lambda} \z \p_\theta(x), \xi_{\theta_n}(x) \right) \dint \mu(x) \right\} \\[1ex]
  &\eqqcolon \tilde{\rho}^{n+1}(\bm{\lambda} \z, \theta),
\end{align*}
where $\tilde{g}(\bm{\lambda} \z) = \gfunc(\z)$ and the induction basis is given by $\rhofunc^0 = \tilde{\rho}^0 = \tilde{g} = \gfunc$.

A necessary condition for $\bm{\lambda}^*$ to solve \eqref{eq:lambda_opt} is that for all $k = 1, \ldots, k$
\begin{equation}
  0 \in \partial_{\lambda_k} \left( \rhofunc(\bm{1}, \theta_0) - \sum_{k=1}^K \lambda_k \overline{\alpha}_k \right) = \partial_{\lambda_k} \rhofunc(\bm{1}, \theta_0) - \overline{\alpha}_k
\end{equation}
By Lemma \ref{lm:diff_lambda_z}, it holds that
\begin{equation}
  \lambda_k \partial_{\lambda_k} \rhofunc(\bm{1}, \theta_0) = \partial_{z_k} \rhofunc(\bm{1},\theta_0)
\end{equation}
so that by definition of $\bm{\lambda}^*$
\begin{equation}
  \lambda_k^* \, \overline{\alpha}_k \in \partial_{z_k} \rho_{\lambda^*}(\bm{1},\theta_0).
\end{equation}
By Theorem \ref{th:rho_derivatives}, it also holds that for all $\pi \in \Pi_{\bm{\lambda}^*}^*(\pmb{\Pbb})$
\begin{equation}
  \lambda_k^* \, \alpha_k(\pi, \Pbb_k) = \lambda_k^* \, \alpha_{\pi,\Pbb_k}^{(k)}(\bm{1}, \theta_0) \in \partial_{z_k} \rho_{\lambda^*}(\bm{1},\theta_0).
\end{equation}
This proves \eqref{eq:constrained_optimal_alpha}. In order to show \eqref{eq:constrained_optimal_gamma}, it suffices to show that $\gamma(\pi^*, \pmb{\Pbb}) \in \partial_{z_k} \rho_{\bm{\lambda}^*}(\bm{1}, \theta_0)$. Lemma~\ref{lm:mixed_policies} guarantees that a policy $\pi^\dagger \in \Pi_{\bm{\lambda}^*}^*(\pmb{\Pbb})$ exists that satisfies the error probability constraints in \eqref{eq:constrained_minimax} with equality, that is
\begin{equation}
  \alpha_k(\pi^\dagger, \Pbb_k) = \overline{\alpha}_k.
\end{equation}
for all $k = 1, \ldots, K$. It then follows that
\begin{align*}
  \gamma(\pi^\dagger, \Pbb_0) + \sum_{k=1}^K \lambda_k^* \overline{\alpha}_k &= \gamma(\pi^\dagger, \Pbb_0) + \sum_{k=1}^K \lambda_k^* \alpha_k(\pi^\dagger, \Pbb_k) \\
  &= \inf_{\pi \in \Pi} \left\{ \gamma(\pi, \Pbb_0) + \sum_{k=1}^K \lambda_k^* \alpha_k(\pi, \Pbb_k) \right\} \\
  &\leq \inf_{\pi \in \PiAlpha} \left\{ \gamma(\pi, \Pbb_0) + \sum_{k=1}^K \lambda_k^* \alpha_k(\pi, \Pbb_k) \right\} \\
  &= \inf_{\pi \in \PiAlpha} \gamma(\pi, \Pbb_0) + \sum_{k=1}^K \lambda_k^* \overline{\alpha}_k,
\end{align*}
which implies $\inf_{\pi \in \PiAlpha} \gamma(\pi, \Pbb_0) = \gamma(\pi^\dagger, \Pbb_0)$. By definition of $\pi^*$ and Theorem~\ref{th:rho_derivatives} it follows that
\begin{equation}
  \gamma(\pi^*, \Pbb_0) = \gamma(\pi^\dagger, \Pbb_0) \in \Big\{ \gamma_{\pi,\Pbb_0}(\bm{1}, \theta_0) : \pi \in \Pi_{\bm{\lambda}^*}^* \Big\} = \partial_{z_0} \rhofunc(\bm{1}, \theta_0).
\end{equation}
This implies that $\pi^\dagger \in \Pi_{\bm{\lambda}^*}^*(\pmb{\Pbb}) \cap \PiAlpha^*(\pmb{\Pbb})$, which concludes the proof.

\section{Proof of Theorem~\titleref{th:lfds_equations}}
\label{apd:proof_lfds_equations}

The existence and uniqueness of $\gamma_{\pi,\Pcal}$ and $\alpha_{\pi,\Pcal}^{(k)}$ can be proven in analogy to the existence and uniqueness of $\rhofunc$ in Theorem~\ref{th:optimal_cost_function}. Consider the sequence of functions $\bigl( \alpha_{\pi,\Pcal}^{(k),n} \bigr)_{n \geq 0}$ that is defined recursively via
\begin{equation}
  \alpha_{\pi,\Pcal}^{(k),n} = \psi(1-\delta_k) + (1-\psi) \left( \sup_{H \in \Pcal_{\theta}} \int \alpha_{\pi,\Pcal}^{(k),n-1} \bigl( \z \p_\theta(x), \xi_\theta(x) \bigr) \, H(\dint x) \right) 
\end{equation}
with $\alpha_{\pi,\Pcal}^{(k),0} = \psi(1-\delta_k)$. It is not hard to show that this sequence is nondecreasing and bounded for all $n \geq 0$. The nondecreasing property can be shown via induction. Assuming $\alpha_{\pi,\Pcal}^{(k),n} \geq  \alpha_{\pi,\Pcal}^{(k),n-1}$, it follows that
\begin{align*}
  \alpha_{\pi,\Pcal}^{(k),n+1} &= \psi(1-\delta_k) + (1-\psi) \left( \sup_{H \in \Pcal_{\theta}} \int \alpha_{\pi,\Pcal}^{(k),n} \bigl( \z \p_\theta(x), \xi_\theta(x) \bigr) \, H(\dint x) \right) \\
  &\geq \psi(1-\delta_k) + (1-\psi) \left( \sup_{H \in \Pcal_{\theta}} \int \alpha_{\pi,\Pcal}^{(k),n-1} \bigl( \z \p_\theta(x), \xi_\theta(x) \bigr) \, H(\dint x) \right) \\
  &= \alpha_{\pi,\Pcal}^{(k),n}.
\end{align*}
The induction basis is given by
\begin{align*}
  \alpha_{\pi,\Pcal}^{(k),1} &=  \psi(1-\delta_k) + (1-\psi) \left( \sup_{H \in \Pcal_{\theta}} \int \alpha_{\pi,\Pcal}^{(k),0} \bigl( \z \p_\theta(x), \xi_\theta(x) \bigr) \, H(\dint x) \right) \alpha_{\pi,\Pcal}^{(k),0} \\
  &\geq \psi(1-\delta_k) = \alpha_{\pi,\Pcal}^{(k),0}
\end{align*}
Boundedness can be shown in the same manner. Assuming that $\alpha_{\pi,\Pcal}^{(k),n} \leq 1$, it follows that
\begin{align*}
  \alpha_{\pi,\Pcal}^{(k),n+1} &= \psi(1-\delta_k) + (1-\psi) \left( \sup_{H \in \Pcal_{\theta}} \int \alpha_{\pi,\Pcal}^{(k),n} \bigl( \z \p_\theta(x), \xi_\theta(x) \bigr) \, H(\dint x) \right) \\
  &\leq \psi(1-\delta_k) + (1-\psi) = 1 - \psi \delta_k \leq 1,
\end{align*}
with induction basis $\alpha_{\pi,\Pcal}^{(k),0} = \psi(1-\delta_k) \leq 1$. Hence, Lemma~\ref{lm:monotonic_convergence} applies and the sequence $\bigl(\alpha_{\pi,\Pcal}^{(k),n}\bigr)_{n \geq 0}$ converges to the unique limit $\alpha_{\pi,\Pcal}^{(k)}$, which satisfies the integral equation \eqref{eq:lfds_equation_errors}.

The same arguments can be used to show existence of $\gamma_{\pi,\Pcal}$, the only difference being that $\gamma_{\pi,\Pcal}$ is not bounded from above. More precisely, the sequence $\bigl( \gamma_{\pi,\Pcal}^n \bigr)_{n \geq 0}$ that is defined recursively via
\begin{equation}
  \gamma_{\pi,\Pcal}^n = (1-\psi) \left( 1 + \sup_{H \in \Pcal_{\theta}} \int \gamma_{\pi,\Pcal}^{n-1}\bigl( \z \p_\theta(x), \xi_\theta(x) \bigr) \, H(\dint x) \right), 
\end{equation}
with $\gamma_{\pi,\Pcal}^n = 0$, can be shown to be nondecreasing. Hence, for every $(\z, \theta) \in \rhodomain$, $\bigl(\gamma_{\pi,\Pcal}^n(\z, \theta)\bigr)_{n \geq 0}$ is a monotonic sequence of real numbers. If this sequence is bounded, the same arguments as before apply and a unique limit $\gamma_{\pi,\Pcal}(\z, \theta)$ exists. If the sequence is unbounded, it is guaranteed to diverge to infinity \cite[Theorem 3.12]{Larson2017}, that is, $\lim_{n \to \infty} \gamma_{\pi,\Pcal}^n(\z, \theta) = \gamma_{\pi,\Pcal}(\z, \theta) = \infty$. Consequently, $\gamma_{\pi,\Pcal}(\z,\theta) \in \Rnng \cup \{ \infty \}$ exists and is unique for every $(\z, \theta) \in \rhodomain$. This concludes the proof.

\section{Proof of Theorem~\titleref{th:lfds}}
\label{apd:proof_lfds}

Theorem~\ref{th:lfds} follows from Theorem~\ref{th:lfds_equations} and can be proven via contradiction. The proof is detailed only for $\Qbb_k$, $k = 1, \ldots, K$; for $k = 0$ it follows analogously. Assume that a distribution $\Pbb^* \in \Pcal_k$ exists such that $\alpha_k(\pi, \Pbb^*) > \alpha_k(\pi, \Qbb_k)$, with $\Qbb_k$ defined in Theorem~\ref{th:lfds}. By \eqref{eq:conditional_unconditional_performace}, this implies that $\alpha_{\pi, \Pbb^*}^{(k)}(\bm{1},\theta_0) > \alpha_{\pi,\Qbb_k}^{(k)}(\bm{1},\theta_0)$, where $\alpha_{\pi, \Pbb^*}^{(k)}$ solves
\begin{equation}
  \alpha_{\pi, \Pbb^*}^{(k)} = \psi (1-\delta_k) + (1-\psi) \int \alpha_{\pi,\Pbb^*}^{(k)} \bigl( \z \p_\theta(x), \xi_\theta(x) \bigr) \, P_\theta^*(\dint x)
\end{equation}
and $\{P_\theta^*\}_{\theta \in \Omega_\theta}$ denotes the family of conditional distributions corresponding to $\Pbb^*$. However, by definition,
\begin{equation}
  \alpha_{\pi, \Pbb^*}^{(k)} \leq \psi (1-\delta_k) + (1-\psi) \sup_{H \in \Pcal_{\theta}^{(k)}} \int \alpha_{\pi,\Pbb^*}^{(k)} \bigl( \z \p_\theta(x), \xi_\theta(x) \bigr) \, H(\dint x)
\end{equation}
so that, using the same arguments as in Appendix~\ref{apd:proof_lfds_equations}, a nondecreasing sequence of functions $\alpha_{\pi,\Pbb^*}^{(k),n}$ can be constructed with $\alpha_{\pi,\Pbb^*}^{(k),0} = \alpha_{\pi,\tilde{\Pbb^*}}^{(k)}$ that converges to $\alpha_{\pi,\Pbb^*}^{(k),\infty}$ for $n \to \infty$. Since by Theorem~\ref{th:lfds_equations} this limit is unique, it follows that $\alpha_{\pi,\Pbb^*}^{(k),\infty} = \alpha_{\pi, \Qbb_k}^{(k)} \geq \alpha_{\pi,\Pbb^*}^{(k)}$, which contradicts the assumption that $\alpha_{\pi,\Pbb^*}^{(k)}(\bm{1},\theta_0) > \alpha_{\pi,\Qbb_k}^{(k)}(\bm{1},\theta_0)$. This concludes the proof.

\section{Proof of Theorem~\titleref{th:minimax_equations}}
\label{apd:proof_minimax_equations}

The proof of Theorem~\ref{th:minimax_equations} closely follows the proof Theorem~5 in \cite{Novikov2009_multiple_hypotheses}. That is, it is shown that the functions $\rhofunc$, $\dfunc$, and $\Dfunc$ can be defined as pointwise limits of monotonic and bounded sequences. From this, existence and uniqueness follow. 

Let $\bigl( \rhofunc^n \bigr)_{n \geq 0}$, $\bigl( \dfunc^n \bigr)_{n \geq 1}$, and $\bigl( \Dfunc^n \bigr)_{n \geq 1}$ be defined recursively via
\begin{align*}
  \rhofunc^n(\z, \theta) &= \min \left\{\, \gfunc(\z) \,,\, z_0 + \dfunc^n(\z, \theta) \,\right\}, \\[1.5ex]
  \dfunc^n(\z, \theta) &= \sup_{\bm{P} \in \bm{\mathcal{P}}_{\theta}} \; \Dfunc^n(\z, \theta; \bm{P}), \\[1ex]
  \Dfunc^n(\z, \theta; \bm{P}) &= \int \rhofunc^{n-1} \bigl( \z \p_\theta(x), \xi_\theta(x) \bigr) \, \mu(\dint x), 
\end{align*}
with $\rhofunc^0 = \gfunc$. Since $\rhofunc^n$ is a nondecreasing function of $\dfunc^n$, $\dfunc^n$ is a nondecreasing function of $\Dfunc^n$, and $\Dfunc^n$ is a nondecreasing function of $\rhofunc^{n-1}$, it follows that all three sequences are nondecreasing. Moreover, since $\rhofunc^n$ is upper bounded by $\gfunc$ for all $n \geq 0$, it holds that
\begin{equation}
  \Dfunc^n(\z, \theta; \bm{P}) \leq \int \gfunc\bigl( \z \p_\theta(x), \xi_\theta(x) \bigr) \, \mu(\dint x) \leq \sum_{k=1}^K \lambda_k z_k
\end{equation}
for all $\bm{P}$ and, consequently,
\begin{equation}
  \dfunc^n(\z, \theta) \leq \sup_{\bm{P} \in \bm{\mathcal{P}}_{\theta}} \; \sum_{k=1}^K \lambda_k z_k = \sum_{k=1}^K \lambda_k z_k.
\end{equation}
This concludes the proof.

\section{Proof of Theorem~\titleref{th:minimax_tests}}
\label{apd:proof_minimax_tests}

A pair $(\pi^*, \pmb{\Qbb})$ is minimax optimal in the sense of \eqref{eq:unconstrained_minimax} if it satisfies the saddle point condition
\begin{equation}
  L_{\bm{\lambda}}(\pi^*, \pmb{\Pbb}) \leq L_{\bm{\lambda}}(\pi^*, \pmb{\Qbb}) \leq L_{\bm{\lambda}}(\pi, \pmb{\Qbb})
  \label{eq:saddle_point_condition}
\end{equation}
for all $\pi \in \Pi$ and all $\pmb{\Pbb} \in \mathbb{M}^{K+1}$. That is, $\pi^*$ is optimal with respect to $\pmb{\Qbb}$, and $\pmb{\Qbb}$ is least favorable with respect to $\pi^*$. 

Assume that the pair $(\pi^*, \pmb{\Qbb})$ satisfies the conditions in Theorem~\ref{th:minimax_tests}. The inequality on the right hand side of \eqref{eq:saddle_point_condition}, that is, optimality of the policy $\pi^*$ w.r.t~$\pmb{\Qbb}$, follows immediately from $\pi^* \in \PiLambda^*(\pmb{\Qbb})$. The inequality on the left hand side, that is, $\pmb{\Qbb}$ being least favorable w.r.t.~$\pi^*$, can be shown as follows. Let $\bm{Q}_{\z,\theta} \in \bm{\Qcal}_{\z,\theta}$ so that
\begin{equation*}
  \Dfunc(\z,\theta; \bm{Q}_{\z,\theta}) = \sup_{\bm{P} \in \bm{\Pcal}_\theta} \Dfunc(\z,\theta; \bm{P}) = \sup_{\bm{P} \in \bm{\Pcal}_\theta} \int \rhofunc \bigl( \z \p(x), \xi_\theta(x) \bigr) \, \mu(\dint x).
\end{equation*}
The partial G\^ateux derivatives \cite{Bell2014} of $\Dfunc$ with respect to $P_k$, $k = 0, \ldots, K$, in the direction $H \in \Pcal_\theta^{(k)}$ are given by
\begin{align*}
  &\lim_{\varepsilon \to 0} \frac{\Dfunc\bigl(\z,\theta; (P_0, \ldots, (1-\varepsilon) P_k + \varepsilon H, \ldots, P_K)\bigr) - \Dfunc(\z,\theta; \bm{P})}{\varepsilon} \\
  &\quad= z_k \int \partial_{z_k} \rhofunc \bigl( \z \p(x), \xi_\theta(x) \bigr) (h(x) - p_k(x)) \, \mu(\dint x) \\
  &\quad = z_k \left( \int \partial_{z_k} \rhofunc \bigl( \z \p(x), \xi_\theta(x) \bigr) \, H(\dint x) - \int \partial_{z_k} \rhofunc \bigl( \z \p(x), \xi_\theta(x) \bigr) \, P_k(\dint x) \right),
\end{align*}
where the limit can be taken inside the integral owing to Leibniz integral rule (Lemma~\ref{lm:leibniz}). A necessary condition for $\bm{Q}_{\z,\theta}$ to be a maximizer of $\Dfunc(\z,\theta; \bm{P})$ is that all partial G\^ateux derivatives evaluated at $\bm{Q}_{\z,\theta}$ are non-positiv in all feasible directions. That is, it holds that
\begin{equation*}
  \int \partial_{z_k} \rhofunc \bigl( \z \q_{\z, \theta}(x), \xi_\theta(x) \bigr) \, Q_{\z,\theta}^{(k)}(\dint x) \geq \int \partial_{z_k} \rhofunc \bigl( \z \q_{\z, \theta}(x), \xi_\theta(x) \bigr) \, H(\dint x)
\end{equation*}
for all $H \in \Pcal_\theta^{(k)}$, or equivalently,
\begin{equation}
  Q_{\z,\theta}^{(k)} \in \argmax_{H \in\Pcal_\theta^{(k)}} \; \int \partial_{z_k} \rhofunc \bigl( \z \q_{\z, \theta}(x), \xi_\theta(x) \bigr) \, H(\dint x).
\end{equation}
By the second statement in Theorem~\ref{th:rho_derivatives}, this implies that for all $\pi \in \Pi_\lambda^*(\pmb{\Qbb})$
\begin{equation}
  Q_{\z,\theta}^{(k)} \in \argmax_{H \in\Pcal_\theta^{(0)}} \; \int \gamma_{\pi,\Pbb_0} \bigl( \z \q_{\z, \theta}(x), \xi_\theta(x) \bigr) \, H(\dint x)
  \label{eq:minimax_lfds_runlength}
\end{equation}
for $k = 0$ and that
\begin{equation}
  Q_{\z,\theta}^{(k)} = \argmax_{H \in\Pcal_\theta^{(k)}} \; \int \alpha_{\pi,\Pbb_k}^{(k)} \bigl( \z \p(x), \xi_\theta(x) \bigr) \, H(\dint x)
  \label{eq:minimax_lfds_errors}
\end{equation}
for $k = 1, \ldots, K$. From \eqref{eq:minimax_lfds_runlength} and \eqref{eq:minimax_lfds_errors} and the Chapman--Kolmogorov equations \eqref{eq:ck_runlength} and \eqref{eq:ck_errors}, it follows that 
\begin{align*}
  \gamma_{\pi^*,\Qbb} &= (1-\psi) \left( 1 + \int \gamma_{\pi,\Qbb}\bigl( \z \q_{\z,\theta}(x), \xi_\theta(x) \bigr) \, Q_{\z,\theta}^{(0)}(\dint x) \right) \\
  &= (1-\psi) \left( 1 + \sup_{H \in \Pcal_{\theta}^{(0)}} \int \gamma_{\pi,\Qbb}\bigl( \z \q_{\z,\theta}(x), \xi_\theta(x) \bigr) \, H(\dint x) \right)
\end{align*}
and
\begin{align*}
  \alpha_{\pi,\Pcal}^{(k)} &= \psi(1-\delta_k) + (1-\psi) \left( \int \alpha_{\pi,\Pcal}^{(k)} \bigl( \z \p_\theta(x), \xi_\theta(x) \bigr) \, Q_{\z,\theta}^{(0)}(\dint x) \right)  \\
  &= \psi(1-\delta_k) + (1-\psi) \left( \sup_{H \in \Pcal_{\theta}^{(k)}} \int \alpha_{\pi,\Pcal}^{(k)} \bigl( \z \p_\theta(x), \xi_\theta(x) \bigr) \, H(\dint x) \right).
\end{align*}
That is, $\gamma_{\pi^*,\Qbb}$ and $\alpha_{\pi^*,\Qbb}^{(k)}$ solve the integral equations \eqref{eq:lfds_equation_runlength} and \eqref{eq:lfds_equation_errors}. By Theorem~\ref{th:lfds}, this implies that $\Qbb_0$ is least favorable with respect to the conditional expected run-length $\gamma_{\pi,\Pbb_0}(\z,\theta)$ for all $(\z, \theta) \in \rhodomain$ and that $\Qbb_k$ is least favorable with respect to the conditional error probabilities $\alpha_{\pi,\Pbb_k}^{(k)}(\z,\theta)$ for all $k = 1, \ldots, K$ and all $(\z, \theta) \in \rhodomain$. Finally, using \eqref{eq:conditional_unconditional_performace}, it follows that
\begin{align*}
  \sup_{\pmb{\Pbb} \in \bm{\Pcal}} L_{\bm{\lambda}}(\pi^*, \pmb{\Pbb}) &= \sup_{\pmb{\Pbb} \in \bm{\Pcal}} \left\{ \gamma(\pi^*,\Pbb_0) + \sum_{k=1}^K \lambda_k \alpha_k(\pi^*,\Pbb_k) \right\} \\
  &= \sup_{\pmb{\Pbb} \in \bm{\Pcal}} \left\{ \gamma_{\pi^*, \Pbb_0}(\bm{1}, \theta_0) + \sum_{k=1}^K \lambda_k \alpha_{\pi^*, \Pbb_k}^{(k)}(\bm{1}, \theta_0) \right\} \\
  &= \gamma_{\pi^*, \Qbb_0}(\bm{1}, \theta_0) + \sum_{k=1}^K \lambda_k \alpha_{\pi^*, \Qbb_k}^{(k)}(\bm{1}, \theta_0) \\
  &= \gamma(\pi^*,\Qbb_0) + \sum_{k=1}^K \lambda_k \alpha_k(\pi^*,\Qbb_k) = L_{\bm{\lambda}}(\pi^*, \pmb{\Qbb}).
\end{align*}
This concludes the proof.

\section{Proof of Theorem~\titleref{th:minimax_constrained_test}}
\label{apd:proof_constrained_minimax_test}

Theorem~\ref{th:minimax_constrained_test} is proven in two steps. First, by definition, all policies $\pi \in \Pi_{\bm{\lambda}^*}^*(\pmb{\Qbb})$ satisfy Theorem~\ref{th:optimal_constrained_test}, which immediately implies \eqref{eq:constrained_minimx_alpha} as well as the existence of a policy $\pi^\dagger \in \Pi_{\bm{\lambda}^*}^*(\pmb{\Qbb})$ that satisfies the constraints on the error probabilities with equality for a given vector of distributions $\pmb{\Qbb}$. The second step is to show that the pair $(\pi^\dagger, \pmb{\Qbb})$ is a saddle point and, hence, minimax optimal. Using the same arguments as in Appendix~\ref{apd:proof_minimax_tests}, it holds that
\begin{align*}
  \gamma(\pi^\dagger, \Qbb_0) + \sum_{k=1}^K \lambda_k^* \overline{\alpha}_k &= \gamma(\pi^\dagger, \Qbb_0) + \sum_{k=1}^K \lambda_k^* \alpha_k(\pi^\dagger, \Qbb_k) \\
  &= \inf_{\pi \in \Pi} \left\{ \gamma(\pi, \Qbb_0) + \sum_{k=1}^K \lambda_k^* \alpha_k(\pi, \Qbb_k) \right\} \\
  &\leq \inf_{\pi \in \PiAlpha} \left\{ \gamma(\pi, \Qbb_0) + \sum_{k=1}^K \lambda_k^* \alpha_k(\pi, \Qbb_k) \right\} \\
  &= \inf_{\pi \in \PiAlpha} \gamma(\pi, \Qbb_0) + \sum_{k=1}^K \lambda_k^* \overline{\alpha}_k \\
  &= \gamma(\pi^*, \Qbb_0) + \sum_{k=1}^K \lambda_k^* \overline{\alpha}_k.
\end{align*}
Moreover, by Theorem~\ref{th:minimax_tests}, it holds that
\begin{align*}
  \gamma(\pi^\dagger, \Qbb_0) + \sum_{k=1}^K \lambda_k^* \overline{\alpha}_k &= \gamma(\pi^\dagger, \Qbb_0) + \sum_{k=1}^K \lambda_k^* \alpha_k(\pi^\dagger, \Qbb_k) \\
  &= \sup_{\pmb{\Pbb} \in \pmb{\Pcal}} \left\{ \gamma(\pi^\dagger, \Pbb_0) + \sum_{k=1}^K \lambda_k^* \alpha_k(\pi^\dagger, \Pbb_k) \right\} \\
  &= \sup_{\Pbb_0 \in \Pcal_0} \; \gamma(\pi^\dagger, \Pbb_0) + \sum_{k=1}^K \lambda_k^* \sup_{\Pbb_k \in \Pcal_k} \; \alpha_k(\pi^\dagger, \Pbb_k) \\
  &= \gamma(\pi^\dagger, \Pbb_0^*) + \sum_{k=1}^K \lambda_k^* \overline{\alpha}_k.
\end{align*}
Hence, $(\pi^\dagger, \pmb{\Qbb})$ satisfies
\begin{equation}
  \inf_{\pi \in \PiAlpha} \; \gamma(\pi, \Qbb_0) = \gamma(\pi^*, \Qbb_0) = \gamma(\pi^\dagger, \Qbb_0) = \gamma(\pi^\dagger, \Pbb_0^*) = \sup_{\Pbb_0 \in \mathbb{M}} \; \gamma(\pi^\dagger, \Pbb_0),
  \label{eq:constrained_minimax_saddlepoint}
\end{equation}
which implies minimax optimality. This concludes the proof.

\section{Proof of Corollary~\titleref{col:constrained_minimax_existence}}
\label{apd:proof_col:constrained_minimax_existence}

Corollary~\ref{col:constrained_minimax_existence} can be proven by a standard duality argument. Let $\PiAlpha(\bm{\Pcal})$ be the set of all policies that satisfy the error probability constraints in \eqref{eq:constrained_minimax}:
\begin{align*}
  \inf_{\pi \in \PiAlpha(\bm{\Pcal})} \; \sup_{\pmb{\Pbb} \in \bm{\Pcal}} \; \gamma(\pi,\Pbb_0) &\geq \sup_{\bm{\lambda} \geq 0} \; \inf_{\pi \in \PiAlpha(\bm{\Pcal})} \; \sup_{\pmb{\Pbb} \in \bm{\Pcal}} \; \left\{ \gamma(\pi,\Pbb_0) + \sum_{k=1}^K \lambda_k \left( \alpha_k(\pi,\Pbb_k) - \overline{\alpha}_k \right) \right\} \\
    &\geq \sup_{\bm{\lambda} \geq 0} \; \left\{ \inf_{\pi \in \Pi} \; \sup_{\pmb{\Pbb} \in \bm{\Pcal}} \; \left\{ \gamma(\pi,\Pbb_0) + \sum_{k=1}^K \lambda_k \alpha_k(\pi,\Pbb_k) \right\} - \sum_{k=1}^K \lambda_k \overline{\alpha}_k \right\} \\
    &= \sup_{\bm{\lambda} \geq 0} \; \left\{ \rhofunc(\bm{1}, \theta_0) - \sum_{k=1}^K \lambda_k \overline{\alpha}_k \right\}.
\end{align*}
Theorem~\ref{col:constrained_minimax_existence} follows immediately.

\bibliographystyle{plain}
\bibliography{bibliography}

\begin{thebibliography}{10}

\bibitem{Medina2015}
M.~Avella~Medina and E.~Ronchetti.
\newblock {R}obust {S}tatistics: {A} {S}elective {O}verview and {N}ew
  {D}irections.
\newblock {\em Wiley Interdisciplinary Reviews: Computational Statistics},
  7(6):372--393, 2015.

\bibitem{Banerjee2015}
T.~Banerjee and V.~V. Veeravalli.
\newblock {D}ata-{E}fficient {M}inimax {Q}uickest {C}hange {D}etection {W}ith
  {C}omposite {P}ost-{C}hange {D}istribution.
\newblock {\em IEEE Transactions on Information Theory}, 61(9):5172--5184,
  2015.

\bibitem{Basu1998}
A.~Basu, I.~R. Harris, N.~L. Hjort, and M.~C. Jones.
\newblock {R}obust and {E}fficient {E}stimation by {M}inimising a {D}ensity
  {P}ower {D}ivergence.
\newblock {\em Biometrika}, 85(3):549--559, 1998.

\bibitem{Bell2014}
J.~Bell.
\newblock {F}r\'{e}chet {D}erivatives and {G}\^{a}teaux {D}erivatives, 2014.
\newblock available online: http://individual.utoronto.ca/jordanbell
  /notes/frechetderivatives.pdf.

\bibitem{Berenguer2015}
I.~M. Berenguer, H.~Kunze, La~D. Torre, and R.~M. Gal{\'a}n.
\newblock {\em {I}nterdisciplinary {T}opics in {A}pplied {M}athematics,
  {M}odeling and {C}omputational {S}cience}, chapter {S}et-{V}alued {N}onlinear
  {F}redholm {I}ntegral {E}quations: {D}irect and {I}nverse {P}roblem, pages
  65--71.
\newblock Springer International Publishing, Basel, Switzerland, 2015.

\bibitem{Breuer2016}
T.~Breuer and I.~Csisz\'{a}r.
\newblock {M}easuring {D}istribution {M}odel {R}isk.
\newblock {\em Mathematical Finance}, 26(2):395--411, 2016.

\bibitem{Brodsky2008a}
B.~E. Brodsky and B.~S. Darkhovsky.
\newblock {M}inimax {M}ethods for {M}ultihypothesis {S}equential {T}esting and
  {C}hange-{P}oint {D}etection {P}roblems.
\newblock {\em Sequential Analysis}, 27(2):141--173, 2008.

\bibitem{Brodsky2008}
B.~E. Brodsky and B.~S. Darkhovsky.
\newblock {M}inimax {S}equential {T}ests for {M}any {C}omposite {H}ypotheses
  {I}.
\newblock {\em Theory of Probability \& Its Applications}, 52(4):565--579,
  2008.

\bibitem{Capasso_Bakstein_2015}
V.~Capasso and D.~Bakstein.
\newblock {\em {A}n {I}ntroduction to {C}ontinuous-{T}ime {S}tochastic
  {P}rocesses}.
\newblock {M}odeling and {S}imulation in {S}cience, {E}ngineering and
  {T}echnology. Birkh{\"a}user Basel, Basel, Switzerland, 3 edition, 2015.

\bibitem{Chieu2009}
N.~H. Chieu.
\newblock {T}he {F}r\'{e}chet and {L}imiting {S}ubdifferentials of {I}ntegral
  {F}unctionals on the {S}paces {$L_1(\Omega, E)$}.
\newblock {\em Journal of Mathematical Analysis and Applications},
  360(2):704--710, 2009.

\bibitem{Choquet1954}
Gustave Choquet.
\newblock {T}heory of {C}apacities.
\newblock {\em Annales de l'Institut Fourier}, 5:131--295, 1954.

\bibitem{Degroot1960}
M.~H. De{G}root.
\newblock {M}inimax {S}equential {T}ests of {S}ome {C}omposite {H}ypotheses.
\newblock {\em The Annals of Mathematical Statistics}, 31(4):1193--1200, 1960.

\bibitem{Dragalin1988}
V.~P. Dragalin and A.~Novikov.
\newblock {A}symptotic {S}olution of the {K}iefer--{W}eiss {P}roblem for
  {P}rocesses with {I}ndependent {I}ncrements.
\newblock {\em Theory of Probability \& Its Applications}, 32(4):617--627,
  1988.

\bibitem{Dvoretzky1953}
{A}. {D}voretzky, {J}. {K}iefer, and {J}. {W}olfowitz.
\newblock {S}equential {D}ecision {P}roblems for {P}rocesses with {C}ontinuous
  {T}ime {P}arameter. {T}esting {H}ypotheses.
\newblock {\em The Annals of Mathematical Statistics}, 24(2):254--264, 1953.

\bibitem{El-Sawy1979}
A.~El-Sawy and V.~D. Vandelinde.
\newblock {R}obust {S}equential {D}etection of {S}ignals in {N}oise.
\newblock {\em IEEE Transactions on Information Theory}, 25(3):346--353, 1979.

\bibitem{Fauss2015}
M.~Fau\ss{} and A.~M. Zoubir.
\newblock {A} {L}inear {P}rogramming {A}pproach to {S}equential {H}ypothesis
  {T}esting.
\newblock {\em Sequential Analysis}, 34(2):235--263, 2015.

\bibitem{Fauss2016_old_bands}
M.~Fau{\ss} and A.~M. Zoubir.
\newblock {O}ld {B}ands, {N}ew {T}racks---{R}evisiting the {B}and {M}odel for
  {R}obust {H}ypothesis {T}esting.
\newblock {\em IEEE Transactions on Signal Processing}, 64(22):5875--5886,
  2016.

\bibitem{Fauss2017_isi}
M.~Fau{\ss} and A.~M. Zoubir.
\newblock {M}inimax {R}obust {S}equential {H}ypothesis {T}esting {U}nder
  {D}ensity {B}and {U}ncertainties.
\newblock In {\em Proc. of the ISI World Statistics Congress (ISI)}, 2017.

\bibitem{Fauss2018_supplement}
M.~Fau\ss{}, A.~M. Zoubir, and H.~V. Poor.
\newblock Supplement to ``{M}inimax {O}ptimal {S}equential {H}ypothesis {T}ests
  for {M}arkov {P}rocesses'', 2018.

\bibitem{Fauss2018}
M.~Fau{\ss}, A.~M. Zoubir, and V.~H. Poor.
\newblock {O}n the {E}quivalence of $f$-{D}ivergence {B}alls and {D}ensity
  {B}ands in {R}obust {D}etection.
\newblock In {\em Proc. of the IEEE International Conference on Acoustics,
  Speech and Signal Processing (ICASSP)}, 2018.

\bibitem{Fauss2016_thesis}
Michael Fau{\ss}.
\newblock {\em {D}esign and {A}nalysis of {O}ptimal and {M}inimax {R}obust
  {S}equential {H}ypthesis {T}ests}.
\newblock PhD thesis, TU Darmstadt, Institute of Communications, 2016.

\bibitem{Fellouris2012}
{G}. {F}ellouris and {A}.~{G}. {T}artakovsky.
\newblock {A}lmost {O}ptimal {S}equential {T}ests of {D}iscrete {C}omposite
  {H}ypotheses, 2012.

\bibitem{Fellouris2012a}
G.~Fellouris and A.~G. Tartakovsky.
\newblock {N}early {M}inimax {O}ne-{S}ided {M}ixture-{B}ased {S}equential
  {T}ests.
\newblock {\em Sequential Analysis}, 31(3):297--325, 2012.

\bibitem{Ferrari2005}
S.~Ferrari and R.~F. Stengel.
\newblock {S}mooth {F}unction {A}pproximation {U}sing {N}eural {N}etworks.
\newblock {\em IEEE Transactions on Neural Networks}, 16(1):24--38, 2005.

\bibitem{Gao2018}
R.~Gao, L.~Xie, Y.~Xie, and H.~Xu.
\newblock {R}obust {H}ypothesis {T}esting {U}sing {W}asserstein {U}ncertainty
  {S}ets.
\newblock {\em Advances in Neural Information Processing Systems}, pages
  7913--7923, 2018.

\bibitem{Ghosh1991}
B.~K. Ghosh and P.~K. Sen, editors.
\newblock {\em {H}andbook of {S}equential {A}nalysis}.
\newblock Statistics: A Series of Textbooks and Monographs. CRC Press, Boca
  Raton, FL, USA, 1991.

\bibitem{Gul2016}
{G}. {G}\"ul and {A}.~{M}. {Z}oubir.
\newblock {R}obust {H}ypothesis {T}esting with $\alpha$-{D}ivergence.
\newblock {\em IEEE Transactions on Signal Processing}, 64(18):4737--4750,
  2016.

\bibitem{GyorfiNemetz1975}
L.~Gy\"orfi and T.~Nemetz.
\newblock {O}n the {D}issimilarity of {P}robability {M}easures.
\newblock Technical report, Mathematical Institute of the Hungarian Academy of
  Science, 1975.

\bibitem{GyorfiNemetz1977}
L.~Gy\"orfi and T.~Nemetz.
\newblock $f$-{D}issimilarity: {A} {G}eneral {C}lass of {S}eparation {M}easures
  of {S}everal {P}robability {D}istributions.
\newblock {\em Colloquia of the J\'anos Bolyai Mathematical Society
  Mathematical Society: Topics in Information Theory}, 16:309--321, 1977.

\bibitem{GyorfiNemetz1978}
L.~Gy\"orfi and T.~Nemetz.
\newblock $f$-{D}issimilarity: {A} {G}eneralization of the {A}ffinity of
  {S}everal {D}istributions.
\newblock {\em Annals of the Institute of Statistical Mathematics},
  30(1):105--113, 1978.

\bibitem{Hafner1993}
R.~Hafner.
\newblock {C}onstruction of {M}inimax-{T}ests for {B}ounded {F}amilies of
  {P}robability-{D}ensities.
\newblock {\em Metrika}, 40(1):1--23, 1993.

\bibitem{EoM1990}
Michiel Hazewinkel, editor.
\newblock {\em {E}ncyclopaedia of {M}athematics}, volume~5, chapter
  {K}olmogorov--{C}hapman {E}quation, page 292.
\newblock Springer, Dordrecht, Netherlands, 1994.

\bibitem{Huber1965}
P.~J. Huber.
\newblock {A} {R}obust {V}ersion of the {P}robability {R}atio {T}est.
\newblock {\em The Annals of Mathematical Statistics}, 36(6):1753--1758, 1965.

\bibitem{Huber1981}
P.~J. Huber.
\newblock {\em {R}obust {S}tatistics}.
\newblock Wiley, Hoboken, NJ, USA, 1981.

\bibitem{Huber1973}
P.~J. Huber and V.~Strassen.
\newblock {M}inimax {T}ests and the {N}eyman--{P}earson {L}emma for
  {C}apacities.
\newblock {\em The Annals of Statistics}, 1(2):251--263, 1973.

\bibitem{Kassam1981}
S.~A. Kassam.
\newblock {R}obust {H}ypothesis {T}esting for {B}ounded {C}lasses of
  {P}robability {D}ensities ({C}orresp.).
\newblock {\em IEEE Transactions on Information Theory}, 27(2):242--247, 1981.

\bibitem{KassamPoor1985}
S.~A. Kassam and H.~V. Poor.
\newblock {R}obust {T}echniques for {S}ignal {P}rocessing: {A} {S}urvey.
\newblock {\em Proceedings of the IEEE}, 73(3):433--481, 1985.

\bibitem{Kharin2002}
A.~Kharin.
\newblock {O}n {R}obustifying of the {S}equential {P}robability {R}atio {T}est
  for a {D}iscrete {M}odel under ``{C}ontaminations''.
\newblock {\em Austrian Journal of Statistics}, 31(4):267--277, 2002.

\bibitem{Kiefer1957}
J.~Kiefer and L.~Weiss.
\newblock {S}ome {P}roperties of {G}eneralized {S}equential {P}robability
  {R}atio {T}ests.
\newblock {\em The Annals of Mathematical Statistics}, 28(1):57--74, 1957.

\bibitem{Kunsch2017}
R.~Kunsch.
\newblock {\em {H}igh-{D}imensional {F}unction {A}pproximation: {B}reaking the
  {C}urse {W}ith {M}onte {C}arlo {M}ethods}.
\newblock PhD thesis, Friedrich-Schiller-Universit{\"a}t Jena, Jena, Germany,
  2017.

\bibitem{Larson2017}
L.~Larson.
\newblock {I}ntroduction to {R}eal {A}nalysis, 2017.
\newblock Lecture Notes.

\bibitem{Lorden1976}
G.~Lorden.
\newblock 2-{SPRT'S} and {T}he {M}odified {K}iefer-{W}eiss {P}roblem of
  {M}inimizing an {E}xpected {S}ample {S}ize.
\newblock {\em The Annals of Statistics}, 4(2):281--291, 1976.

\bibitem{Maronna2006}
R.~Maronna, D.~Martin, and V.~Yohai.
\newblock {\em {R}obust {S}tatistics: {T}heory and {M}ethods}.
\newblock Wiley, Hoboken, NJ, USA, 2006.

\bibitem{Maurice1957}
R.~J. Maurice.
\newblock {A} {M}inimax {P}rocedure for {C}hoosing {B}etween {T}wo
  {P}opulations using {S}equential {S}ampling.
\newblock {\em Journal of the Royal Statistical Society. Series B
  (Methodological)}, 19(2):255--261, 1957.

\bibitem{Mukherjee2005}
M.~N. Mukherjee.
\newblock {\em {E}lements of {M}etric {S}paces}.
\newblock Academic Publishers, Kolkata, India, 2005.

\bibitem{Nguyen2009}
X.~Nguyen, M.~J. Wainwright, and M.~I. Jordan.
\newblock {O}n {S}urrogate {L}oss {F}unctions and $f$-{D}ivergences.
\newblock {\em The Annals of Statistics}, 37(2):876--904, 2009.

\bibitem{Novikov2009_multiple_hypotheses}
A.~Novikov.
\newblock {O}ptimal {S}equential {M}ultiple {H}ypothesis {T}ests.
\newblock {\em Kybernetika}, 45(2):309--330, 2009.

\bibitem{Oesterreicher1978}
F.~\"{O}sterreicher.
\newblock {O}n the {C}onstruction of {L}east {F}avourable {P}airs of
  {D}istributions.
\newblock {\em Zeitschrift f\"ur Wahrscheinlichkeitstheorie und Verwandte
  Gebiete}, 43(1):49--55, 1978.

\bibitem{Papageorgiou1997}
N.~S. Papageorgiou.
\newblock {C}onvex {I}ntegral {F}unctionals.
\newblock {\em Transactions of the American Mathematical Society},
  349(4):1421--1436, 1997.

\bibitem{Pardo2005}
L.~Pardo.
\newblock {\em {S}tatistical {I}nference {B}ased on {D}ivergence {M}easures}.
\newblock CRC Press, Boca Raton, FL, USA, 2005.

\bibitem{Pavlov1991}
I.~V. Pavlov.
\newblock {S}equential {P}rocedure of {T}esting {C}omposite {H}ypotheses with
  {A}pplications to the {K}iefer--{W}eiss {P}roblem.
\newblock {\em Theory of Probability \& Its Applications}, 35(2):280--292,
  1991.

\bibitem{Poor1980}
H.~V. Poor.
\newblock {R}obust {D}ecision {D}esign {U}sing a {D}istance {C}riterion.
\newblock {\em IEEE Transactions on Information Theory}, 26(5):575--587, 1980.

\bibitem{Poor2009}
H.~V. Poor and O.~Hadlijiadis.
\newblock {\em {Q}uickest {D}etection}.
\newblock Cambridge University Press., Cambridge, UK, 2009.

\bibitem{Reid2011}
M.~D. Reid and R.~C. Williamson.
\newblock {I}nformation, {D}ivergence and {R}isk for {B}inary {E}xperiments.
\newblock {\em Journal of Machine Learning Research}, 12:731--817, 2011.

\bibitem{Rockafellar1968}
R.~T. Rockafellar.
\newblock {I}ntegrals {W}hich are {C}onvex {F}unctionals.
\newblock {\em Pacific Journal of Mathematics}, 24(3):525--539, 1968.

\bibitem{Rockafellar1970}
R.~T. Rockafellar.
\newblock {\em {C}onvex {A}nalysis}.
\newblock Princeton University Press, Princeton, NJ, USA, 1970.

\bibitem{Rockafellar1974}
R.~T. Rockafellar.
\newblock {\em {C}onjugate {D}uality and {O}ptimization}, chapter~1, pages
  1--74.
\newblock Society for Industrial and Applied Mathematics, Philadelphia, PA,
  USA, 1974.

\bibitem{Rudin1976}
W.~Rudin.
\newblock {\em {P}rinciples of {M}athematical {A}nalysis}.
\newblock McGraw-Hill, New York City, NY, USA, 1976.

\bibitem{Schmitz1987}
N.~Schmitz.
\newblock {M}inimax {S}equential {T}ests of {C}omposite {H}ypotheses on the
  {D}rift of a {W}iener {P}rocess.
\newblock {\em Statistische Hefte}, 28(1):247--261, 1987.

\bibitem{Siegmund1985}
D.~Siegmund.
\newblock {\em {S}equential {A}nalysis}.
\newblock Springer, New York City, NY, USA, 1985.

\bibitem{Sochman2005}
J.~Sochman and J.~Matas.
\newblock {W}ald{B}oost---{L}earning for {T}ime {C}onstrained {S}equential
  {D}etection,.
\newblock In C.~Schmid, S.~Soatto, and C.~Tomasi, editors, {\em Proc. of the
  IEEE Computer Society Conference on Computer Vision and Pattern Recognition},
  pages 150--156, 2005.

\bibitem{Tartakovsky2014}
A.~Tartakovsky, I.~Nikiforov, and M.~Basseville.
\newblock {\em {S}equential {A}nalysis: {H}ypothesis {T}esting and
  {C}hangepoint {D}etection}.
\newblock Chapman and Hall/CRC, Boca Raton, FL, USA, 2014.

\bibitem{Tartakovsky2003}
A.~G. Tartakovsky, X.~R. Li, and G.~Yaralov.
\newblock {S}equential {D}etection of {T}argets in {M}ultichannel {S}ystems.
\newblock {\em IEEE Transactions on Information Theory}, 49:425--445, 2003.

\bibitem{Unnikrishnan2011}
J.~Unnikrishnan, V.~V. Veeravalli, and S.~P. Meyn.
\newblock {M}inimax {R}obust {Q}uickest {C}hange {D}etection.
\newblock {\em IEEE Transactions on Information Theory}, 57(3):1604--1614,
  2011.

\bibitem{Varshney2011}
K.~R. Varshney.
\newblock {B}ayes {R}isk {E}rror is a {B}regman {D}ivergence.
\newblock {\em IEEE Transactions on Signal Processing}, 59(9):4470--4472, 2011.

\bibitem{Verdu1984}
S.~Verd\'u and H.~Poor.
\newblock {O}n {M}inimax {R}obustness: {A} {G}eneral {A}pproach and
  {A}pplications.
\newblock {\em IEEE Transactions on Information Theory}, 30(2):328--340, 1984.

\bibitem{Vos2001}
H.~J. Vos.
\newblock {A} {M}inimax {P}rocedure in the {C}ontext of {S}equential {T}esting
  {P}roblems in {P}sychodiagnostics.
\newblock {\em British Journal of Mathematical and Statistical Psychology},
  54:139--159, 2001.

\bibitem{Voudouri1988}
E.~Voudouri and L.~Kurz.
\newblock {A} {R}obust {A}pproach to {S}equential {D}etection.
\newblock {\em IEEE Transactions on Acoustics, Speech and Signal Processing},
  36(8):1200--1210, 1988.

\bibitem{Wald1947}
A.~Wald.
\newblock {\em {S}equential {A}nalysis}.
\newblock Wiley, Hoboken, NJ, USA, 1947.

\bibitem{Weiss2013}
C.~H. Wei\ss{} and H.-Y. Kim.
\newblock {B}inomial {AR}{(1)} {P}rocesses: {M}oments, {C}umulants, and
  {E}stimation.
\newblock {\em Statistics}, 47(3):494--510, 2013.

\bibitem{Yilmaz2016}
E.~Y{\i}lmaz, J.~F. Gemmeke, and H.~Van Hamme.
\newblock {N}oise {R}obust {E}xemplar {M}atching {W}ith {A}lpha–{B}eta
  {D}ivergence.
\newblock {\em Speech Communication}, 76:127--142, 2016.

\bibitem{Zhitlukhin2013}
M.~V. Zhitlukhin, A.~A. Muravlev, and A.~N. Shiryaev.
\newblock {T}he {O}ptimal {D}ecision {R}ule in the {K}iefer--{W}eiss {P}roblem
  for a {B}rownian {M}otion.
\newblock {\em Russian Mathematical Surveys}, 68(2):389--391, 2013.

\bibitem{Zoubir2012}
A.~M. Zoubir, V.~Koivunen, Y.~Chakhchoukh, and M.~Muma.
\newblock {R}obust {E}stimation in {S}ignal {P}rocessing: {A}
  {T}utorial-{S}tyle {T}reatment of {F}undamental {C}oncepts.
\newblock {\em IEEE Signal Processing Magazine}, 29(4):61--80, 2012.

\bibitem{Zoubir2018}
A.~M. Zoubir, V.~Koivunen, E.~Ollila, and M.~Muma.
\newblock {\em {R}obust {S}tatistics for {S}ignal {P}rocessing}.
\newblock Cambridge University Press, Cambridge, UK, 2018.

\end{thebibliography}

\end{document}